\documentclass[12pt,a4paper]{article}
\usepackage[T1]{fontenc}
\usepackage{amsmath}
\usepackage{amssymb}
\usepackage{graphicx}
\usepackage[titletoc]{appendix}
\usepackage{hyperref}
\usepackage{xypic}
\usepackage{appendix}

\newcommand{\mZ}{{\mathbb Z}}
\newcommand{\AG}{{\mathcal{ABG}}}
\newcommand{\BM}{{\mathcal{BM}}}

\newcommand{\angA}{\mathcal{B}(A)}
\newcommand{\angAD}{\mathcal{B}(A(D))}
\newcommand{\innAD}{{\text{Inner}(A(D))}}
\newcommand{\angab}{\langle a,b\rangle}
\newcommand{\angcd}{{\langle c,d\rangle}}

\newcommand{\LieR}{{\text{\bf Lie}_{\mathfrak R}}}
\newcommand{\Alt}{{\text{\bf Alt}_k}}

\newcommand{\BAD}{\mathcal{B}(A(D))}
\newcommand{\AGAD}{\mathcal{ABG}(A(D))}
\newcommand{\AGADo}{\mathcal{ABG}(A(1+D))}
\newcommand{\Hom}{{\text{Hom}}}

\newtheorem{cor}{\noindent\rm\bf Corollary}[section]
\newtheorem{thm}{\noindent\rm\bf Theorem}[section]

\newtheorem{pro}{\noindent\rm\bf Proposition}[section]
\newtheorem{defn}{\noindent\rm\bf Definition}[section]

\newtheorem{lem}{\noindent\rm\bf Lemma}[section]
\newtheorem{conj}{\noindent\rm\bf Conjecture}
\newtheorem{conje}{\noindent\rm\bf Conjecture}

\begin{document}
	
\title{Allison-Benkart-Gao functor and the cyclicity of free alternative functors}
\author{Shikui Shang}
\date{}

\maketitle

\begin{center}{Department of Mathematics,	Shanghai University,\\ Shanghai 200444,	PR.China.\\
E\_mail:skshang@mail.ustc.edu.cn}\end{center}

\vspace{2mm}

\abstract{Let $k$ be a field of characteristic $0$. We introduce a pair of adjoint functors, Allison-Benkart-Gao functor $\AG$ and Berman-Moody functor $\BM$, between the category of non-unital alternative algebras over $k$ and the category $\LieR$ of Lie algebras with compatible $sl_3(k)$-actions. Surprisingly, when $A$ is an alternative algebra without a unit, the Allison-Benkart-Gao Lie algebra $\AG(A)$ is not isomorphic to the more well-known Steinberg Lie algebra $st_3(A)$ in general.
	
Let $A(D)$ be the free (non-unital) alternative algebra over $D$ generators with the inner derivation algebra $\innAD$. A conjecture on the homology $H_r(\AGAD)$ is proposed. Furthermore,
consider the degree $n$ component of  $A(D)_n$(resp. $\innAD_n$). The previous conjecture implies another conjecture on the dimensions on $A(D)_n$ and $\text{Inner} A(D)_n$. Some evidences are given to support these conjectures.
 
Finally, we prove the cyclicity of the alternative structure, namely that the
symmetric group $S_{1+D}$ acts on the multilinear part of $A(D)$, which plays an important role to connect the Lie algebra homology of
$\AGAD$ and the character of $A(D)$.
}

{\bf 2020 Mathematics Subject Classification:}primary 17D05; secondary 17B55, 17B60. 

{\bf Keywords:} Free alternative algebra; Allison-Benkart-Gao functor;  the cyclicity of the alternative functor.

\tableofcontents

\section{Introduction}

Let  $k$ be a field of characteristic $0$.  The tensor products $\otimes$ are taken over $k$ if there are no subscripts in this paper.  

Let $A$ be a non-associative algebra over $k$. The associator of $a,b,c\in A$ is defined by
$$(a,b,c)=(ab)c-a(bc).$$

$A$ is called an \emph{alternative algebra} over $k$ if
$$(a,a,b)=0,\ \ \ \ (a,b,b)=0,$$
hold for any $a,b\in A$.

The theory of alternative algebras is one of important parts of the theory of non-associative algebras. Historically, some significant progresses were achieved by M.Zorn \cite{Z}, Schafer \cite{S1}, Skornyakov \cite{Sk}, Shirshov \cite{Shi} and E. Kleinfeld \cite{Kl1,Kl2} and so on. Systematic theories on alternative algebras can be found in the two classical books \cite{S2} and \cite{ZSSS}. Some topics on free alternative algebras were studied by  A. V. Il'tyakov \cite{I} and I. P. Shestakov \cite{She1,She2}, et al.

Recently, some relations between Jordan algebras and Lie algebra of $sl_2$-type are explored by I. Kashuba and O. Mathieu in \cite{KM}. A series of conjectures is given to develop a new viewpoint on (non-unital) free Jordan algebras from the homology of TAG Lie algebras associated with them as coordinate algebras. Specifically, the  dimensions of the $n$-homogeneous components of free Jordan algebras (and their inner derivation algebras) are encoded as the exponents in some infinite product equations. See the various versions of the Conjecture 1 in \cite{KM}.

It is well known that alternative algebras can be viewed as the coordinate algebras of Lie algebras of $sl_3$-type($A_2$-type). Their relations can be found in \cite{ABG}-\cite{BM}. For any alternative algebra $A$, the Lie algebra $\mathcal{ABG}(A)$
is introduced by the functorial construction in \cite{ABG} and \cite{BGKN}. We call it Allison-Benkart-Gao functor $\mathcal{ABG}$ following the works of these authors. The functor $\mathcal{ABG}$ and its adjoint functor $\mathcal{BM}$, Berman-Moody functor, are studied detailedly in this paper.

Next, we focus on the free alternative algebras $A(D)$ with $D$ generators. Some conjectures are proposed, involving the homology of $\AGAD$, their classes of $GL(D)$-modules in the Grothendieck ring and the dimensions of the homogeneous components of $A(D)$ and its inner derivation algebra $\innAD$, which are analogous to those conjectures raised in \cite{KM}. 

An analytic functor is called cyclic if it is the image of another functor under the suspension. In \cite{KM}, the functors $J$ and $T$ of free Jordan and free associative algebras are proved to be cyclic. It is also true for the free alternative functor $A$. This means that there exist actions of symmetric group $S_{1+D}$ on the multilinear part of $A(D)$ and it is another example of "hidden $\Sigma_{n+1}$-actions" first shown by O.Mathieu in \cite{M}. Applying the cyclicity, we can relate $A(D)$ to $\AGADo$ and obtain a conjecture declaring that the $sl_3(k)$-invariant components of $H_r(\AGAD)$ are vanishing for $r\geq 1$, which implies the previous conjecture on the equivalence classes of the Grothendick ring.

The paper is organized as follows. In section 2, we give the details of the construction of the adjoint pair, Allison-Benkart-Gao functor and Berman-Moody functor. Section 3 investigates the homologies $H_r(\AGAD)$ of Lie algebra $\AGAD$ when $0\leq r\leq 2$ and a conjecture(Conjecture \ref{con:1}) on the higher $H_r(\AGAD)$ is stated. In section 4, the classes of $A(D)$ and $\innAD$ in the Grothendieck ring ${\mathcal R}(GL_n(D))$ are considered and two related conjectures(Conjecture \ref{con:2},\ref{con:3}) are given, the relations among which are studied. In section 5, some numerical and non-numerical evidences are shown to support these conjectures. In the main part of Section 6, we prove that $A$ is a cyclic functor and the fact that all $H_r(\AGADo)^{sl_3}=0(r\geq 1)$ is already sufficient to imply conjecture \ref{con:2}.  Finally, we summarize the main results and give some comments in the last section of this paper.

\section{The Allison-Benkart-Gao functor and Berman-Moody functor}
\subsection{The Allison-Benkart-Gao construction}
Let $A$ be an alternative algebra over $k$ which is not necessarily unital.  Denote 
$L_a$ and $R_a$ are respectively the left and right multiplication on $A$ by $a\in A$. A $k$-linear morphism $D\in\text{End}_k(A)$ is called a \emph{derivation}
on $A$ if for all $a,b\in A$,
$D(ab)=D(a)b+aD(b)$,  which is equivalent to $[D,L_a]=L_{D(a)}$ or $[D,R_b]=R_{D(b)}$ in $\text{End}_k(A)$.

For $a,b\in A$, define
$$D_{a,b}=[L_a,L_b]+[R_a,R_b]+[L_a,R_b]\in\text{End}_k(A).$$

The following results are well-known (can be found in \cite{BGKN} or \cite{S2})
\begin{pro}\label{pro:2.1}  For $a,b,c\in A$, $D\in\text{Der}(A)$,

(1) $D_{a,b}=L_{[a,b]}-R_{[a,b]}-3[L_a,R_b]$, $D_{a,b}(c)=[[a,b],c]+3(a,c,b)$;

(2) $D_{a,b}\in\text{Der}(A)$.

(3) $D_{a,b}+D_{b,a}=0$;

(4) $D_{ab,c}+D_{bc,a}+D_{ca,b}=0$;

(5) $[D,D_{a,b}]=D_{D(a),b}+D_{a,D(b)}$.
\end{pro}

 A linear span of $D_{a,b}$ is called an \emph{inner derivation} of $A$, and we denote the subspace of all inner derivations by  $\text{Inner}(A)$, which is an ideal of $\text{Der}(A)$ by (5) of the previous proposition.

\ 

Let $A$ be an alternative algebra. Set $I(A)$ is the subspace of $\Lambda^2A$ spanned by
$$ab\wedge c+bc\wedge a+ca\wedge b,
$$
where $a,b,c$ run over the elements in $A$. Set the quotient space $$\angA=\Lambda^2A/I(A),$$
and denote the canonical image of $a\wedge b$ in $\angA$ by $\angab$.

Let $$sl_3(k)=\{x\in M_{3\times 3}(k)\ |\ \text{Tr}(x)=0\}$$
be the split simple Lie algebra over $k$ of type $A_2$. The Cartan subalgebra ${\frak h}$  is the set of all the traceless diagonal matrices in $sl_3(k)$.

\

{\bf Allison-Benkart-Gao Construction} 

In \cite{ABG}, from the appropriate unital coordinate algebras according to the type of root systems, a unified explicit construction of centrally closed Lie algebras graded by finite root systems is expressed systematically by B. Allison, G. Benkart and Y. Gao, which is functorial on the coordinate algebras. We apply this construction for non-unital alternative algebras to construct Lie algbras of $A_2$-type and call it Allison-Benkart-Gao functor in this paper.

For an alternative algebra $A$,
the Allison-Benkart-Gao algebra $\AG(A)$ is defined as the vector space
$$\AG(A)=\angA\oplus(sl_3(k)\otimes A).$$

The bracket product is defined by bilinear
extension of the formulae
\begin{align*}
	1.\ \ \ \ [x\otimes a,y\otimes b]=&\frac{\text{Tr}(xy)\langle a,b\rangle}{3}+[x,y]\otimes\frac{ab+ba}{2}\\
	&+\left(xy+yx-\frac{2}{3}\text{Tr}(xy)I_3\right)\otimes\frac{ab-ba}{2},\\	
	2.\ \ \ \ [\angab, x\otimes c]=&-[x\otimes c,\angab]=x\otimes D_{a,b}(c),\\
	3.\ \ \ \ [\angab,\angcd]=&\langle D_{a,b}(c),d\rangle+\langle c,D_{a,b}(d)\rangle,
\end{align*}
where $\text{Tr}$ is the trace form on $sl_3(k)$ and $a,b,c,d\in A, x,y\in sl_3(k)$ . 

In the Theorem 4.13 in \cite{ABG} or Proposition 2.2 in \cite{BGKN}, one has seen that
\begin{thm}\label{thm:2.2} $\AG(A)$ is a Lie algebra over $k$ for any alternative algebra $A$.  $\Box$
\end{thm}

{\bf Remark. }When $A$ is a unital alternative algebra, it is proved that $\AG(A)$ is isormorphic to the Steinberg Lie algebra $st_3(A)$ in \cite{BGKN}. We will see that is not true when $A$ is non-unital in the observation after Lemma \ref{lem:3.1}.

\subsection{The Berman-Moody functor}

\

We introduce the category $\LieR$ in this subsection. 

Let $\text{R}$ be the subcategory of Lie algebras over $k$, whose objects are also $sl_3(k)$-modules and the elements of $sl_3(k)$ act as derivations on them. The homomorphisms in this category are morphisms of both Lie algebras and $sl_3(k)$-modules simultaneously. Note that any Lie algebra containing a subalgebra isomorphic to $sl_3(k)$ is an object in $\text{R}$.

Let $$\Delta=\{\alpha_1,\alpha_2,\alpha_1+\alpha_2,-\alpha_1,-\alpha_2,-(\alpha_1+\alpha_2)\}$$
be the root system of $sl_3(k)$. Let $\mathfrak R$ be the subcategory of $sl_3(k)$-modules whose objects are weight $sl_3(k)$-modules with weights only in 
$\Delta\cup\{0\}$.  If $M$ is an object in $\mathfrak R$, then $M=M^{sl_3}\oplus M^{\text{ad}}$, where $M^{sl_3}$ is the invariant submodule of $M$ and $M^{\text{ad}}$ is the isotypic component of $M$ of adjoint type.

\begin{defn}\label{def:2.1} 
	Let $\LieR$ be the full subcategory of $\text{R}$ such that all objects in $\LieR$
	are also objects in the category $\mathfrak R$, namely that any object $\frak g$ in $\LieR$ has a decomposition
	$${\mathfrak g}={\mathfrak g}^{sl_3}\oplus{\mathfrak g}^{\text{ad}},$$ 
	as an $sl_3(k)$-module, 
\end{defn}
Moreover, the weight space decomposition of $\mathfrak g$ is given by
$${\frak g}={\frak g}_0\oplus\bigoplus_{\alpha\in\Delta}{\frak g}_\alpha.$$

In \cite{BM}, the Lie algebras graded by finite root systems are introduced. All $A_2$-graded Lie algebras are objects in $\LieR$, which contain split simple subalgebras isomorphic to $sl_3(k)$.   

\begin{lem}\label{lem:2.1}  If $A$ is an alternative algebra over $k$, then the Allison-Benkart-Gao Lie algebra $\AG(A)$ is an object in $\LieR$. And, $\AG:\Alt\rightarrow\LieR, A\mapsto\AG(A)$ is a functor from the category $\Alt$ of alternative algebras over $k$ to the category $\LieR$. Indeed, one has that $$\AG(A)^{sl_3}=\angA,\ \ \ \AG(A)^{\text{ad}}=sl_3(k)\otimes A.$$
\end{lem}
{\bf Proof. }The actions of $sl_3(k)$ are given by
$$x.(y\otimes a)=[x,y]\otimes a,\ \ \ x.\angab=0,$$
for $x,y\in sl_3(k)$ and $a,b\in A$.

It is straightforward to show that the actions are derivations for the brackets given in the definition of $\AG(A)$.  

The fact that
$A\mapsto\angA$ is functorial implies that $\AG(-)={\mathcal B}(-)\oplus(sl_3(k)\otimes-)$ is also a functor.
$\Box$ 

Note that if $A$ is not unital, $\AG(A)$ is not an $A_2$-graded Lie algebra since the split simple subalgebra does not always exist. Hence, the category $\LieR$ contains more objects besides $A_2$-graded Lie algebras. At least, some Lie algebras in $\LieR$ are not perfect.

\ 

For a Lie algebra $\mathfrak g$ in $\LieR$, denote 
$$\BM({\mathfrak g})={\frak g}_{\alpha_1+\alpha_2},$$
as the highest weight space in ${\mathfrak g}$. 
Then,
$${\mathfrak g}={\mathfrak g}^{sl_3}\oplus(sl_3(k)\otimes \BM({\frak g})).$$

When $\frak g$ is an $A_2$-graded Lie algebras, the Recognition Theorem for types $A_2$ in Section 3 of \cite{BM} shows that $\BM({\frak g})$ is an alternative algebra with identity over $k$ such that 
$$[e_{ij}\otimes a,e_{kl}\otimes b]=\delta_{jk}e_{il}\otimes ab-\delta_{li}e_{kj}\otimes ba$$ 
for $1\leq i,j,k,l\leq 3, i\neq j,k\neq l,\{i,j\}\neq\{k,l\}$. 

Respectively, in category $\LieR$, one has
\begin{thm}\label{thm:2.3} For $\frak g\in\text{Obj}(\LieR)$, $\BM({\mathfrak g})$ is an alternative algebra over $k$ satisfying for $	1\leq i,j,k,l\leq 3, i\neq j,k\neq l,\{i,j\}\neq\{k,l\}$, $a,b\in\BM({\frak g})$,
	\begin{equation}[e_{ij}\otimes a,e_{kl}\otimes b]=\delta_{jk}e_{il}\otimes ab-\delta_{li}e_{kj}\otimes ba.\tag{2.1}\end{equation}
	Furthermore, $\BM$ is a functor from the category $\LieR$ to the category $\Alt$.
\end{thm}
{\bf Proof. }First, we see that
the actions $e_i,f_i(i=1,2)\in sl_3(k)$ on $\frak g$ are nilpotent. One can replace the adjoint actions of the split simple subalgebra in an $A_2$-graded Lie algebra by the derivation actions of $sl_3(k)$ and obtain that there is a product on $\BM({\frak g})$ such that the equation $(2.1)$
holds.

Next, we will show the product is alternative.
The proof of Recognition Theorem in \cite{BM} has to be modified on account of the absence of the split simple subalgebra . We consider the actions of $sl_3(k)$ instead of the existence of the identity in $\BM(\frak g)$. 

For $a,b,c\in\BM({\frak g})$, Jacobi identity gives that
\begin{align*}
	&[[e_{13}\otimes a,e_{31}\otimes b], e_{23}\otimes c]=[e_{13}\otimes a,[e_{31}\otimes b, e_{23}\otimes c]]\\
	=&-[e_{13}\otimes a,e_{21}\otimes cb]=e_{23}\otimes(cb)a.\tag{2.2}
\end{align*}
Then, by the derivation action of $e_{12}$ and $(2.2)$,
\begin{align*}
	&[[e_{13}\otimes a,e_{31}\otimes b], e_{13}\otimes c]=[[e_{13}\otimes a,e_{31}\otimes b], e_{12}.(e_{23}\otimes c)]\\
	=&e_{12}.e_{23}\otimes(cb)a-[e_{12}.([e_{13}\otimes a,e_{31}\otimes b]), e_{23}\otimes c]\\
	=&e_{13}\otimes(cb)a+[[e_{13}\otimes a,e_{32}\otimes b],e_{23}\otimes c]\\
	=&e_{13}\otimes((ab)c+(cb)a).\tag{2.3}
\end{align*}
On the other hand, by Jacobi identity,
\begin{align*}
	&[[e_{13}\otimes a,e_{31}\otimes b], e_{12}\otimes c]=[e_{13}\otimes a,[e_{31}\otimes b, e_{12}\otimes c]]\\=&e_{12}\otimes a(bc).\tag{2.4}
\end{align*}
And, by the derivation action of $e_{23}$ and $(2.4)$
\begin{align*}
	&[[e_{13}\otimes a,e_{31}\otimes b], e_{13}\otimes c]=-[[e_{13}\otimes a,e_{31}\otimes b], e_{23}.(e_{12}\otimes c)]\\
	=&-e_{23}.e_{12}\otimes a(bc)+[e_{23}.([e_{13}\otimes a,e_{31}\otimes b]), e_{12}\otimes c]\\
	=&e_{13}\otimes a(bc)+[[e_{13}\otimes a,e_{21}\otimes b],e_{12}\otimes c]\\
	=&e_{13}\otimes(a(bc)+c(ba)).\tag{2.5}
\end{align*}
Hence, $(ab)c+(cb)a=a(bc)+c(ba)$ by $(2.3)$ and $(2.5)$, i.e.,
\begin{equation*}(a,b,c)+(c,b,a)=0.\tag{2.6}\end{equation*}
Meanwhile, by $[e_{13}\otimes a,e_{12}\otimes b]=0, [e_{13}\otimes a,e_{23}\otimes c]=0$ and the derivation action of $e_{31}$,
\begin{align*}
	&[e_{31}.e_{13}\otimes a, e_{13}\otimes bc]=[e_{31}.e_{13}\otimes a, [e_{12}\otimes b,e_{23}\otimes c]]\\
	=&[[e_{31}.e_{13}\otimes a,e_{12}\otimes b],e_{23}\otimes c]+[e_{12}\otimes b,[e_{31}.e_{13}\otimes a, e_{23}\otimes c]]\\
	=&-[[e_{13}\otimes a,e_{31}.e_{12}\otimes b],e_{23}\otimes c]-[e_{12}\otimes b,[e_{13}\otimes a, e_{31}.e_{23}\otimes c]]\\
	=&-e_{13}\otimes((ab)c+b(ca)).\tag{2.7}
\end{align*}
And, by $e_{12}.e_{13}\otimes a=0$ and $[e_{13}\otimes a,e_{23}\otimes bc]=0$,
\begin{align*}
	&[e_{31}.e_{13}\otimes a, e_{13}\otimes bc]=[e_{31}.e_{13}\otimes a, e_{12}.e_{23}\otimes bc]\\
	=&e_{12}.[e_{31}.e_{13}\otimes a, e_{23}\otimes bc]-[e_{12}.e_{31}.e_{13}\otimes a, e_{23}\otimes bc]\\
	=&-e_{12}.[e_{13}\otimes a, e_{31}.e_{23}\otimes bc]+[e_{32}.e_{13}\otimes a, e_{23}\otimes bc]\\
	=&-e_{13}((bc)a+a(bc)).\tag{2.8}
\end{align*}
We also obtain that $(ab)c+b(ca)=(bc)a+a(bc)$ by $(2.7)$ and $(2.8)$, i.e.,
\begin{equation*}(a,b,c)=(b,c,a).\tag{2.9}\end{equation*}

By $(2.6)$ and $(2.9)$, we have that $A$ is an alternative algebra over $k$. 

Finally, if $\phi:{\frak g}_1\rightarrow{\frak g}_2$ is a homomorphism in $\LieR$, $\phi$ is a morphism of $sl_3(k)$-modules. $\phi$ preserves the weight spaces. Particularly, restricting $\phi$ on the highest weight space
gives a $k$-linear morphism from $\BM({\frak g}_1)$ to $\BM({\frak g}_2)$. It is a morphism of alternative algebras since $\phi$ is also a morphism of Lie algebras. This prove the funtoriality of $\BM$.
$\Box$

\

Next, we follow the steps in \cite{BGKN} to show that there is a homomorphism from $\AG(\BM({\frak g}))$ to $\frak g$ in $\LieR$. The trouble is also the absence of an identity in $\BM({\frak g})$. To overcome this obstacle, the actions of $sl_3(k)$ are indispensable.

Note that ${\frak g}={\frak g}^{sl_3}\oplus sl_3\otimes{\BM(\frak g)}$.
For $a,b\in\BM(\frak g),1\leq i\neq j\leq 3$, let
$$H_{ij}(a,b)=[e_{ij}\otimes a, e_{ji}\otimes b].$$
Then, for $k\neq i,j$,
\begin{align*}H_{ij}(a,b)=&[e_{ij}\otimes a, e_{ji}\otimes b]=[e_{ij}\otimes a, e_{jk}.e_{ki}\otimes b]\\
	=&e_{jk}.[e_{ij}\otimes a, e_{ki}\otimes b]-[e_{jk}.e_{ij}\otimes a, e_{ki}\otimes b]\\
	=&H_{ik}(a,b)-(e_{jj}-e_{kk})\otimes ba.\end{align*}
Hence, taking $i=1$, \begin{align*}&H_{1j}(a,b)-(e_{11}-e_{jj})\otimes ba\\=&H_{1k}(a,b)-(e_{jj}-e_{kk})\otimes ba-(e_{11}-e_{jj})\otimes ba\\
	=&H_{1k}(a,b)-(e_{11}-e_{kk})\otimes ba\end{align*}

For $j\in\{2,3\}$, define
$$h(a,b)=H_{1j}(a,b)-(e_{11}-e_{jj})\otimes ba,$$
which does not depend on the choices of $j=2$ or $3$. And,

\begin{lem}\label{lem:2.3} For $a,b,c\in\BM(\frak g)$, the associator $(a,b,c)$,
	\begin{align*}(1)\ \ \ \ \ &h(a,b)+h(b,a)=0,\\
		(2)\ \ \ \ \ &h(ab,c)+h(bc,a)+h(ca,b)=\\
		&(e_{11}-e_{22})\otimes(a,b,c)+(e_{11}-e_{33})\otimes(a,b,c).\end{align*}
\end{lem}
{\bf Proof.} (1) Using the derivations in $sl_3(k)$,
\begin{align*}&h(a,b)=H_{12}(a,b)-(e_{11}-e_{22})\otimes ba\\
	=&[e_{13}.e_{32}\otimes a,e_{21}\otimes b]-(e_{11}-e_{22})\otimes ba\\
	=&e_{13}.e_{31}\otimes ab-[e_{32}\otimes a,e_{13}.e_{21}\otimes b]-(e_{11}-e_{22})\otimes ba\\
	=&(e_{11}-e_{33})\otimes ab-H_{23}(b,a)-(e_{11}-e_{22})\otimes ba.\tag{2.10}\end{align*}
Meanwhile,
\begin{align*}&h(b,a)=H_{13}(b,a)-(e_{11}-e_{33})\otimes ab\\
	=&[e_{12}.e_{23}\otimes b,e_{31}\otimes a]-(e_{11}-e_{33})\otimes ab\\
	=&(e_{11}-e_{22})\otimes ba+H_{23}(b,a)-(e_{11}-e_{33})\otimes ab.\tag{ 2.11}\end{align*}
Hence, $h(a,b)+h(b,a)=0$ from $(2.10)$ and $(2.11)$.

(2) Follow the proof of Lemma 2.14 in \cite{BGKN} and replace any symbol $H_{ij}(1,a)$ by $(e_{ii}-e_{jj})\otimes a\in{\frak g}$, for $1\leq i\neq j\leq 3$. We get
\begin{align*}&h(ab,c)-h(a,bc)+h(ca,b)=\\
	&(e_{11}-e_{22})\otimes(a,b,c)+(e_{11}-e_{33})\otimes(a,b,c).\end{align*}
Applying (1), we have the identity (2). $\Box$ 

We now define a bilinear form on $\BM({\frak g})$
$$\psi:\BM({\mathfrak g})\times\BM({\mathfrak g})\rightarrow{\mathfrak g}$$ such that
\begin{align*}\psi(a,b)=&3h(a,b)-(e_{11}-e_{22})\otimes(ab-ba)\\
	&-(e_{11}-e_{33})\otimes(ab-ba).\end{align*}

\begin{lem}\label{lem:2.4} The image of $\psi$ is contained in ${\frak g}^{sl_3}$. And,
	\begin{align*}(1)\ \ \ \ \ &\psi(a,b)+\psi(b,a)=0,\\
		(2)\ \ \ \ \ &\psi(ab,c)+\psi(bc,a)+\psi(ca,b)=0.\end{align*}
\end{lem}
{\bf Proof.} For the first assertion, it is enough to show that $$e_{12}.\psi(a,b)=0, \ \ \ \ e_{23}.\psi(a,b)=0.$$
In fact,
\begin{align*}&e_{12}.\psi(a,b)=3e_{12}.([e_{13}\otimes a,e_{31}\otimes b]-(e_{11}-e_{33})\otimes ba)\\
	-&[e_{12},e_{11}-e_{12}]\otimes(ab-ba)-[e_{12},e_{11}-e_{33}]\otimes(ab-ba)\\
	=&-3e_{12}\otimes ab+3e_{12}\otimes ba+2e_{12}\otimes(ab-ba)+e_{12}\otimes(ab-ba)\\
	=&e_{12}\otimes(-3ab+3ba+2ab-2ba+ab-ba)=0.\tag{2.12}\end{align*}
And,
\begin{align*}&e_{23}.\psi(a,b)=3e_{23}.([e_{13}\otimes a,e_{31}\otimes b]-(e_{11}-e_{33})\otimes ba)\\
	-&[e_{23},e_{11}-e_{12}]\otimes(ab-ba)-[e_{23},e_{11}-e_{33}]\otimes(ab-ba)\\
	=&e_{23}\otimes(-3ba+3ba-ab+ba+ab-ba)=0.\tag{2.13}\end{align*}

Next, also modifying the proof of Lemma 2.16 in \cite{BGKN}, one can get (1) and (2). $\Box$ 

Hence, $\psi$ induces a linear map $\tilde{\psi}$ on the quotient ${\mathcal B}(\BM({\mathfrak g}))$, $$\tilde{\psi}:{\mathcal B}(\BM({\mathfrak g}))\rightarrow{\mathfrak g}^{sl_3}.$$ And,
one can check that
\begin{align*}[x\otimes a,y\otimes b]=&\frac{1}{3}\text{Tr}(xy)\tilde{\psi}( a,b)+[x,y]\otimes\frac{ab+ba}{2}\\
	&+\left(xy+yx-\frac{2}{3}\text{Tr}(xy)I_3\right)\otimes\frac{ab-ba}{2}\end{align*}
for $x,y\in sl_3(k)$ and $a,b\in \BM({\mathfrak g})$.

\begin {lem}\label{lem:2.5}Let ${\mathfrak g}\in\text{Obj }(\LieR)$. Then, there is a homomorphism in $\LieR$
$$\theta_{\mathfrak g}:\AG(\BM({\mathfrak g}))\rightarrow{\mathfrak g}$$
which is the identity restricted on $\BM(\mathfrak g)$.
\end {lem}

{\bf Proof. }The linear map $\theta_{\mathfrak g}:\AG(\BM({\mathfrak g}))\rightarrow{\mathfrak g}$ is defined by requiring that $\theta_{\mathfrak g}$ is identity on $sl_3\otimes \BM({\mathfrak g})$ and $\theta_{\mathfrak g}=\tilde{\psi}$ on ${\mathcal B}(\BM({\mathfrak g}))$. We have that $\theta_{\mathfrak g}$
is a homomorphism in the category $\LieR$.   $\Box$ 

\

The following theorem gives the relation between the functor $\AG$ and $\BM$.
\begin{thm}\label{thm:2.4} The functor $\AG:\Alt\rightarrow\LieR$ is the left adjoint of the Berman-Moody functor $\BM:\LieR\rightarrow\Alt$, namely, there exists a natural equivalence
	$$\Hom_\LieR(\AG(A),{\mathfrak g})\simeq_N\Hom_\Alt(A,\BM({\mathfrak g}))$$
	for any $A\in\Alt$ and ${\mathfrak g}\in\LieR$.
\end{thm}
{\bf Proof. } Since $\BM(\AG(A))=A$, there is a natural map by restricting on $A$
$$\mu:\Hom_{\LieR}(\AG(A), {\mathfrak g})\rightarrow\Hom_{\Alt}(A, \BM({\mathfrak g})), \eta\mapsto\eta|_A.$$
Since $\AG(A)$ is generated by $sl_3\otimes A$, $\mu$ is injective. 

On the other hand, let $\phi:A\rightarrow\BM(\mathfrak g)$
be a morphism of alternative algebras. By functoriality of the $\AG$-construction, we get a morphism of Lie algebras
$$\AG(\phi):\AG(A)\rightarrow\AG(\BM({\mathfrak g}))$$
and there is a canonical Lie algebra morphism by Lemma \ref{lem:2.5}
$$\theta_{\mathfrak g}:\AG(\BM({\mathfrak g}))\rightarrow{\mathfrak g}.$$

So, $\theta_{\mathfrak g}\circ\AG(\phi)$ extends $\phi$ to a morphism of Lie algebras. Therefore, $\phi=\mu(\theta_{\mathfrak g}\circ\AG(\phi))$ and $\mu$ is bijective.
$\Box$

\section{The homology of Lie algebra $\AG(A(D))$}
\subsection{The (co)homology of Lie algebras in $\LieR$}

Let $\mathfrak g$ be a Lie algebra over $k$. 
For $r\geq 0$, $r$-th homologies (resp. cohomologies) $H_r({\mathfrak g},-)$(resp. $H^r({\mathfrak g},-)$) are functors from the category of $\mathfrak g$-modules to the category of $k$-spaces. $H_r({\mathfrak g},M)$(resp. $H^r({\mathfrak g},M)$) is called the $r$-th homology(resp. cohomology) space of $\mathfrak g$ with coefficients in $M$. In particular, we call $H_r({\mathfrak g})=H_r({\mathfrak g},k)$(resp. $H^r({\mathfrak g})=H^r({\mathfrak g},k)$) the $r$-th homology(resp. cohomology) space of $\mathfrak g$ for the trivial $\mathfrak g$-module $k$.

The (co)homology of ${\mathfrak g}$ can be computed using the Chevalley-Eilenberg complex $(\wedge^*{\frak g}, d_*)$, where the differentials $d_n:\wedge^{r+1}{\frak g}\rightarrow\wedge^r{\frak g}$ given by
\begin{equation*}d_r(x_0\wedge\cdots\wedge x_r)=\sum_{i=0}^{r-1}(-1)^ix_0\wedge\cdots\wedge[x_i,x_{i+1}]\wedge\cdots\wedge x_r.\tag{3.1}\end{equation*}
More details can be found in Chapter 7 of \cite {W}.

\

We introduce the suitable module category for Lie algebras in $\LieR$.
\begin {defn}\label{def:3.1}
Let ${\mathfrak g}$ be a Lie algebra in the category $\text{R}$. A ${\mathfrak g}$-module $M$ is called admissible, if it is also an $sl_3(k)$-module and these actions are compatible, i.e., for any $a\in{\frak g}, x\in sl_3(k)$ and $m\in M$, one has
\begin{equation*}
(x.a).m=x.(a.m)-a.(x.m). \tag{3.2}
\end{equation*}
holds. All admissible ${\mathfrak g}$-modules make up a category ${\mathcal M}({\frak g})$ of modules, whose homomorphisms are the $k$-linear morphisms commuting with both $\frak g$-actions and $sl_3(k)$-actions simultaneously. Moreover, the category of admissible $\mathfrak g$-modules in $\mathfrak R$ is denoted by ${\mathcal M}_{\mathfrak R}({\frak g})$ 
\end{defn}

There are examples of admissible $\frak g$-modules,

(1) The direct sums, the tensor products and the quotient of  admissible modules are admissible; 

(2) The adjoint module ${\frak g}$ is admissible. Then, the tensor product $T({\frak g})=\oplus_{n\geq 0}T^n({\frak g})$ and the universal enveloping algebra $U(\mathfrak g)$ are admissible;

(3) If $M$ is an $sl_3(k)$-module, we define the trivial actions of $\frak g$ on $M$. Then, $M$ is admissible. In particular, the trivial module $L(0)$ is admissible.

Using above constructions, a free $\mathfrak g$-module generated by $sl_3(k)$-modules is admissible. Furthermore, any $\mathfrak g$-module generated by $sl_3(k)$-modules subject to relations stable under the $sl_3(k)$-actions is always admissible.

\

Consider the (co)homology of Lie algebras in $\LieR$ with the coefficients of admissible $\frak g$-modules,
\begin {pro}\label{pro:3.2}
(1) For $r\geq 0$, all $H_r({\mathfrak g},-)$ and $H^r({\mathfrak g},-)$
are functors from the category ${\mathcal M}_{\mathfrak R}({\mathfrak g})$ to the category of  $sl_3(k)$-modules. 

(2) If $M$ is a weight $sl_3(k)$-module, so are $H_r({\mathfrak g},M)$ and $H^r({\mathfrak g},M)$. Moreover, if the highest weight of $M$ exists and is $\mu$,  the highest weight of $H_r({\mathfrak g},M)$ is at most
$r(\alpha_1+\alpha_2)+\mu$.
\end {pro}
{\bf Proof. }(1) Consider the Chevalley-Eilenberg complex. It is well-known that the tensor products and wedge products of $sl_3(k)$-modules are also $sl_3(k)$-modules. 
Hence, so are $\wedge^r{\frak g}\otimes M$ for all $r\geq 0$. 

Next, since the actions of $sl_3(k)$  are derivations on $\frak g$ and $M$ is admissible, all the differentials $d_r$ in $(\wedge^*{\frak g}\otimes M,d_*)$ commute with the actions of $sl_3(k)$ by $(3.1)$ and $(3.2)$. That is, $d_r$ are morphisms of $sl_3(k)$-modules. Hence, the kernels and images of $d_r$ are $sl_3(k)$-modules and so are $H_r({\frak g},M)$. The discussion for $H^r({\frak g},M)$ are similar.

(2) is straightforward.  $\Box$

Since the trivial module $k$ is admissible, both $H_r({\mathfrak g})$ and $H^r({\mathfrak g})$ are weight modules of $sl_3(k)$ and their highest weights are at most $r(\alpha_1+\alpha_2)$ for $r\geq 0$.

\

For Lie algebra $\AG(A)$, one has 
\begin{lem}\label{lem:3.1} For any alternative algebra $A$, as an $sl_3(k)$-module
$$H_1({\AG(A)})\simeq L(\alpha_1+\alpha_2)\otimes(A/A^2),$$
which yields that $H_1({\AG(A)})$ has only isotopic component of adjoint type.
\end{lem}
{\bf Proof.} Using the Chevalley-Eilenberg complex, we have $H_1({\mathfrak g})={\mathfrak g}/[{\mathfrak g},{\mathfrak g}]$. For ${\mathfrak g}=\AG(A)$, 
it is easy to see that
$$[\AG(A),\AG(A)]=\angA\oplus(sl_3(k)\otimes A^2).$$
Hence, $$H_1({\AG(A)})=sl_3(k)\otimes(A/A^2)\simeq L(\alpha_1+\alpha_2)\otimes(A/A^2)$$ as an $sl_3(k)$-module.    $\Box$

{\bf A surprising observation} It is not hard to show that for Steinberg Lie algebra $st_3(A)$,
$$st_3(A)={\frak T}\oplus\bigoplus_{1\leq i\neq j\leq 3}X_{ij}(A),$$
where ${\frak T}=\sum_{1\leq i\neq j\leq 3}[X_{ij}(A),X_{ji}(A)]$.
Hence,
$$H_1(st_3(A))=\oplus_{1\leq i\neq j\leq 3}X_{ij}(A/A^2)$$
which is not an $sl_3(k)$-module when $A/A^2\neq0$. Hence, $st_3(A)$ is not an object in $\LieR$ for the reason that there are no appropriate $sl_3(k)$-actions and
it can not be isomorphic to $\AG(A)$.

\subsection{The homology of $\AG(A(D))$}

We consider the free non-unital alternative algebra $A(D)$ with $D$ generators $\{x_1,\cdots, x_D\}$ of in this subsection.

The second (co)homology space of a Lie algebra is related to the central extensions of it.

Let ${\frak g}$ be a Lie algebra in category $\LieR$. Consider the central extension of Lie algebras
$$0\rightarrow M\rightarrow\hat{\frak g}\rightarrow{\frak g}\rightarrow0$$
with a $2$-cocycle $\omega:{\frak g}\times{\frak g}\rightarrow M$. $\hat{\frak g}$ is contained in $\LieR$ if it has compatible $sl_3(k)$-module and $\frak g$-module structures. It is equivalent to that $M$ is admissible and $\omega$ is $sl_3(k)$-invariant, i.e.,
\begin{align*}
	&(x.a).m=x.(a.m)-a.(x.m),\\
	&x.\omega(a,b)=\omega(x.a,b)+(a,x.b)
\end{align*}
for all $x\in sl_3(k)$ and $a,b\in{\frak g}, m\in M$. That is, $\omega\in Z^2({\frak g},M)^{sl_3}$.

\begin {lem}\label{lem:3.3} Let $M$ be a finite-dimensional admissible $\AGAD$-module. Then,
$$H_2(\AG(A(D)),M)^{sl_3}=0.$$
\end {lem}
{\bf Proof.} By duality, we only need to show $H^2(AG(A(D)),M^*)^{sl_3}=0$.

For any $sl_3(k)$-invariant $2$-cocycle $\omega:\AGAD\times\AGAD\rightarrow M^*$,
let $L=M^*\oplus\AGAD$ be the Lie algebra associated to $\omega$.
We have that $L$ is contained in $\LieR$. By Theorem \ref {thm:2.4}, the exact sequence of Lie algebras in category $\LieR$
$$0\rightarrow M^*\rightarrow L\rightarrow\AGAD\rightarrow0$$
splits. Therefore, $H^2(\AGAD,M^*)^{sl_3}=0$. $\Box$ 

\

Here, we summarize the information of $H_r(\AGAD)$ for $r=0,1,2$,
\begin {pro}\label{pro:3.1}(1) $H_0(\AGAD)\simeq L(0)$,

(2) $H_1(\AGAD)\simeq L(\alpha_1+\alpha_2)\otimes k^D$, 

(3) $H_2(\AGAD)^{sl_3}=0,\ \ H_2(\AGAD)^{ad}=0$.

\end {pro}
{\bf Proof. }(1)It is true for any Lie algebra over $k$. 

(2)Obviously, $A(D)^2=\oplus_{n\geq 2}A(D)_n$ and $$A(D)/A(D)^2\simeq A(D)_1.$$ Note that
$A(D)_1\simeq k^D$, where $k^D$ is the $k$-space with dimension $D$. By Lemma \ref{lem:3.1},
we obtain (2).

(3) Taking $M$ is the trivial module $k$  in Lemma \ref{lem:3.3}, we have that  $H_2(\AGAD)^{sl_3}=0$.

Furthermore, let $L(\alpha_1+\alpha_2)$ be the adjoint $sl_3(k)$-module.
With the trivial action of $\AGAD$, $L(\alpha_1+\alpha_2)$ is admissible. Also by above Lemma,
$$H_2(\AGAD,L(\alpha_1+\alpha_2))^{sl_3}=0.$$

Since the $\AGAD$-actions are trivial, we have that $$H_*(\AGAD,L(\alpha_1+\alpha_2))\simeq H_*(\AGAD)\otimes L(\alpha_1+\alpha_2).$$ 
If $H_2(\AGAD)^{\text{ad}}\neq 0$, it has a submodule $N$ isomorphic to $L(\alpha_1+\alpha_2)$ and
$N\otimes L(\alpha_1+\alpha_2)$ is a submodule of $H_2(\AGAD)\otimes L(\alpha_1+\alpha_2)$. However, by the decomposition
\begin{align*}&L(\alpha_1+\alpha_2)^{\otimes 2}\simeq\\
	& L(2\alpha_1+2\alpha_2)\oplus L(2\alpha_1+\alpha_2)\oplus L(\alpha_1+2\alpha_2)\oplus2L(\alpha_1+\alpha_2)\oplus L(0),\end{align*}
$H_2(\AGAD,L(\alpha_1+\alpha_2))$ has a submodule isomorphic to $L(0)$, which conflicts the fact that $H_2(\AGAD,L(\alpha_1+\alpha_2))^{sl_3}=0$. Therefore, we have that $H_2(\AGAD)^{\text{ad}}=0$.
$\Box$  

\

It is a well-known fact that if $\frak m$ is an ordinary free Lie algebra, then $H_r({\frak m})=0$ for $r\geq 2$(See Corollary 7.2.5 in \cite{W}). Here, $A(D)$ are free objects in the category $\Alt$ and $\AGAD$ are free objects in the category $\LieR$ by Theorem \ref{thm:2.4}. Following the viewpoint in \cite{KM}, since only the trivial and adjoint modules of $sl_3$-type occur in Lie algebras of $\LieR$, we give a natural conjecture on the homology modules of the Lie algebra $\AGAD$.
\begin {conj}\label{con:1}As $sl_3(k)$-modules,
$$H_r(\AGAD)^{sl_3}=0 \text{\ \ \  and   }\ \ \ H_r(\AGAD)^{\text{ad}}=0$$
for any $r\geq 2$.
\end {conj}

We have seen that the conjecture holds when $r=2$ by Proposition \ref{pro:3.1}(3). One also has
\begin {lem}\label{lem:3.4}Assume that $H_k(\AG(A(D)))^{sl_3}=0$ for any odd $k$. Then, we have that
$$\angAD=\text{Inner}\ A(D).$$
\end {lem}
{\bf Proof. }We can follow the proof of Lemma 12 in \cite{KM}, where they apply a long exact sequence of cohomology spaces. The version for cohomology of groups can be found in \cite{Ho}. $\Box$

\section{The homogenous dimensions of $A(D)$ and $\angAD$}

Throughout this section, we focus on the dimensions of the homogeneous spaces of $A(D)$ and $\angAD$. Fixing $D\geq 1$, there exist $\mZ_{>0}$-gradings
$$A(D)=\oplus_{n>0}A(D)_n,\ \ \ \angAD=\oplus_{n>0}\angAD_n.$$

For $n\geq 1$, set
$$a_n(D)=\dim A(D)_n, \ \ b_n(D)=\dim\angAD_n.$$
What have we known about $a_n(D)$?
Straightforwardly, $a_1(D)=D$. By Artin's theorem, $a_2(D)=D^2$. We will see that $a_3(D)=D^3+{D\choose 3}$ in Section 5.3. Furthermore, $a_n(D)$ and $b_n(D)$ are closely related to the homology of the Lie algebra $\AGAD$.

\

\subsection{The Grothendieck rings and $\lambda$-operations}

First, we consider the general linear group $G=GL(D)$ over $k$. The analytic $G$-modules are described in Section 1 of \cite{KM} or Chapter 1 of \cite{Ma}.  

Let $Z=k^*\text{id}$ be the subgroup of $G$ constituting by nonzero scale matrices and $V$ an analytic $G$-module. For a non-negative integer $n$, set
$$V_n=\{v\in V\ |\ h.v=\lambda^nv,\forall h\in Z\},$$
for $h=\lambda\cdot\text{id}\in Z$. Then,
$V=\oplus_{n\geq 0}V_n$ and $V_n$ is the homogenous components in $V$ of degree $n$. Indeed, $A(D)$, $\angAD$ and $\innAD$ are analytic $G$-modules.

\

Let $K$ be a reductive algebraic subgroup of $G$, which contains $Z$.
For $n\geq 0$, a rational $K$-module is called a $K$-module of degree $n$, if $\lambda\cdot\text{id}$ acts by $\lambda^n$ for any $\lambda\in k^*$. Denote $\text{Rep}_n(K)$ be the category of $K$-modules of degree $n$.
Let $K_0(\text{Rep}_n(K))$ be the Grothendieck group of the category of homogeneous polynomial $G$-module with degree $n$. Set
\begin{align*}
	&{\mathcal R}(K)=\sum_{n\geq 0}K_0(\text{Rep}_n(K)),\\
	&{\mathcal M}_{>m}(K)=\sum_{n>m}K_0(\text{Rep}_n(K)),\\
	&{\mathcal M}(K)={\mathcal M}_{>0}(K).
\end{align*}

Note that there are products $$K_0(\text{Rep}_m(K))\times K_0(\text{Rep}_n(K))\rightarrow K_0(\text{Rep}_{m+n}(K)),$$
induced by the tensor product of the $K$-modules.
${\mathcal R}(K)$ becomes a commutative ring and ${\mathcal M}(K)$ is an ideal.
Moreover, ${\mathcal R}(K)$ is complete with respect to the ${\mathcal M}(K)$-adic topology. Any element $a$ of ${\mathcal R}(K)$ can be written as a formal series $a =\sum_{n\geq 0}a_n$ where $a_n\in K_0(\text{Rep}_n(K))$. 

Particularly, taking $K=Z$, there is an isomorphism of rings $${\mathcal R}(Z)\rightarrow \mZ[[z]], [V]\mapsto\sum_{n\geq 0}\dim(V_n)z^n,$$ with ${\mathcal M}(Z)$ identifying $z\mZ[[z]]$ and ${\mathcal R}(Z)/{\mathcal M}(Z)\simeq\mZ$. In the rest of this paper, we always identify the ring ${\mathcal R}(Z)$ with $\mZ[[z]]$.

In \cite{KM}, we have seen that there exists a $\lambda$-operation on ${\mathcal M}(K)$, which is a homomorphism of groups from the additive group ${\mathcal M}(K)$ to the multiplicative group $1+{\mathcal M}(K)$. That is, for any $a_1,a_2\in{\mathcal M}(K)$, one has  $\lambda(a_1+a_2)=\lambda(a_1)\lambda(a_2)$. And, for the effective class $a=[V]$ in ${\mathcal M}(K)$,
$$\lambda(a)=\sum_{r\geq 0}(-1)^r[\Lambda^rV]\in1+{\mathcal M}(K).$$

\subsection{$A(D)$ and $\angAD$ in the Grothendieck rings}

Let $P_+$ be the subset of $\mZ^2_{\geq0}$ with dominant integral weights in the root lattice of type $A_2$. Then,
$$P_+=\{(m_1,m_2)\in\mZ^2_{\geq0}\ |\  2m_1-m_2\geq 0,2m_2-m_1\geq 0\}.$$
Let $L(m_1\alpha_1+m_2\alpha_2)$ be the finite-dimensional irreducible $PSL(3)$-module of highest weight $m_1\alpha_1+m_2\alpha_2$ for $(m_1,m_2)\in P_+$. Then, $$\{[L(m_1\alpha_1+m_2\alpha_2)]\ | (m_1,m_2)\in P_+\}$$ is a $\mZ$-basis of $K_0(PSL(3))$. Since $K_0(K\times PSL(3))=K_0(K)\otimes_\mZ K_0(PSL(3))$, any element $a\in K_0(K\times PSL(3))$ can be written as a sum
$$a=\sum_{(m_1,m_2)\in P_+}[a:L(m_1\alpha_1+m_2\alpha_2)][L(m_1\alpha_1+m_2\alpha_2)],$$
where $[a:L(m_1\alpha_1+m_2\alpha_2)]\in{\mathcal R}(K)$. Particularly, $L(0)$ is the the trivial module and $L(\alpha_1+\alpha_2)$ is the adjoint module of $sl_3(k)$. And,
$$[\AGAD]=[\BAD][L(0)]+[A(D)][L(\alpha_1+\alpha_2])]$$

The following Lemma is critical for our purpose.
\begin {lem}\label{lem:4.1} (Main Lemma for alternative case) Let $K$ be an algebraic subgroup of $G$ which contains $Z$.

(1) There are unique $a(K),b(K)\in{\mathcal M}(K)$ such that the element $$a(K)[L(\alpha_1+\alpha_2)]+b(K)\in{\mathcal M}(K\times PSL(3))$$ satisfies the following two equations :
\begin{align*}
	&[\lambda\left(a(K)[L(\alpha_1+\alpha_2)]+b(K)\right):L(0)]=[k|_K],\\
	&[\lambda\left(a(K)[L(\alpha_1+\alpha_2)]+b(K)\right):L(\alpha_1+\alpha_2)]=-[k^D|_K],
\end{align*}
where $k$ is the trivial $G$-module, $k^D$ is the natural $G$-module and $k|_K$ and $k^D|_K$ are their restrictions on the subgroup $K$.

(2)For $K=G$, set ${\mathcal A}(D)=a(G)$ and ${\mathcal B}(D)=b(G)$ in ${\mathcal M}(G)$. 

Then, $a(K)={\mathcal A}(D)|_K$ and $b(G)={\mathcal B}(D)|_K$.
\end {lem}
{\bf Proof. }(1) This follows immediately from the analogous assertion on the Jordan algebra case, which is a proof of Hensel's Lemma type in a complete ring. See the proof of Lemma 1 in \cite{KM}.

(2) Use the uniqueness of  $a(K)$ and $b(K)$.   $\Box$

\

Then, we can give the main conjecture for the structure of free alternative algebra $\AG(A(D))$.
\begin {conj}\label{con:2} Let $D$ be a non-negative integer and $A(D)$ the free alternative algebra with $D$ generators. In ${\mathcal R}(G)$, we have that
$${\mathcal A}(D)=[A(D)],\ \ {\mathcal B}(D)=[\angAD],$$
where ${\mathcal A}(D)$ and ${\mathcal B}(D)$ are defined in the above lemma (2).
\end {conj}

This conjecture is not surprise by the following Theorem.
\begin{thm}\label{thm:4.1} Conjecture \ref{con:1} implies Conjecture \ref{con:2}.
\end{thm}
{\bf Proof.} First, note that in ${\mathcal R}(K\times PSL(3))$
$$[\AGAD]=[\angAD][L(0)]+[A(D)][L(\alpha_1+\alpha_2)].$$
By the definition of $\lambda$-operation,
$$\lambda([\AGAD])=\sum_{r\geq 0}(-1)^r[\wedge^r(\AGAD)].$$

Meanwhile, using Euler's characteristic formula for  Chevalley-Eilenberg complex of $\AGAD$, we have
$$\lambda([\AGAD])=\sum_{r\geq 0}(-1)^r[H_r(\AGAD)].$$

Assume that Conjecture \ref{con:1} holds, by Proposition \ref{pro:3.1},
\begin{align*}
	&\lambda([\AGAD])^{sl_3}=[H_0(\AGAD)^{sl_3}]=[k],\\
	&\lambda([\AGAD])^{\text{ad}}=-[H_1(\AGAD)^{\text{ad}}]=-[L(\alpha_1+\alpha_2)\otimes k^D].
\end{align*}
We have that
\begin{align*}
	&[\lambda([\AGAD]):L(0)]=[k],\\
	&[\lambda([\AGAD]):L(\alpha_1+\alpha_2)]=-[k^D].
\end{align*}

By the uniqueness of ${\mathcal A}(D)$ and ${\mathcal B}(D)$ in Lemma \ref{lem:4.1}, we get
$${\mathcal A}(D)=[A(D)],\ \ {\mathcal B}(D)=[\angAD].\ \ \Box$$

\subsection{The dimensions of $A(D)_n$ and $\angAD_n$}
In this section, we take $K=Z$ and give a weak form of conjecture 2 on the dimension of $A(D)$  and $\angAD$.

We have identified the ring ${\mathcal R}(Z)$ with $\mZ[[z]]$. And, $K_0(PSL(3))$ corresponds the subring $\mZ[t_1^{\pm 1},t_2^{\pm 1}]^{S_3}$ in $\mZ[t_1^{\pm 1},t_2^{\pm 1}]$ consisting of the symmetric polynomials under the symmetric group $S_3$, as the Weyl group of the $A_2$-type, where $t_1=e^{\alpha_1},t_2=e^{\alpha_2}$ for the simple roots $\alpha_1,\alpha_2$ of the simple algebraic group $PSL(3)$. It follows that
$${\mathcal R}(Z\times PLS(3))=\mZ[t_1^{\pm 1},t_2^{\pm 1}]^{S_3}[[z]].$$

For $n>1$, set the dimensions
$$\dim A(D)_n=a_n(D), \ \ \dim \angAD_n=b_n(D),$$
and the two power series $$a(z)=\sum_{n\geq 1}a_n(D)z^n, \ \ \ b(z)=\sum_{n\geq 1}b_n(D)z^n\in z\mZ[[z]].$$

If Conjecture \ref{con:2} holds, using Lemma \ref{lem:4.1}(2) for $K=Z$, $a(z)$ and $b(z)$ are uniquely determined by
the two equations in $\mZ[[z]]$:
\begin{align*}
	&[\lambda\left(a(z)[L(\alpha_1+\alpha_2)]+b(z)\right):L(0)]=1,\\
	&[\lambda\left(a(z)[L(\alpha_1+\alpha_2)]+b(z)\right):L(\alpha_1+\alpha_2)]=-Dz.\tag{4.1}
\end{align*}

\

For two series $a(z)=\sum_{n\geq 1}a_nz^n$, $b(z)=\sum_{n\geq 1}b_nz^n$ in $z\mZ[[z]]$, define the power series of the infinite product form
\begin{align*}
	&\Phi(z)=\Phi(a(z),b(z))\\
	=&\prod_{n\geq 1}\left(\prod_{\alpha\in\Delta}(1-e^\alpha z^n)^{a_n}\right)(1-z^n)^{2a_n+b_n}\\
	&\in1+z\mZ[t_1^{\pm 1},t_2^{\pm 1}]^{S_3}[[z]].\tag{4.2}
\end{align*}

By Proposition A of 22.5 in \cite{Hu},
we can set $$\Phi(a(z),b(z))=\sum_{m,n\geq 0}\Phi(m_1\alpha_1+m_2\alpha_2)\text{ch}\ L(m_1\alpha_1+m_2\alpha_2),$$
where $\Phi(m_1\alpha_1+m_2\alpha_2)\in\mZ[[z]]$ can be viewed as the multiplicity of $L(m_1\alpha_1+m_2\alpha_2)$.

\begin{lem}\label{lem:4.4}There are uniquely two power series $a(z)$ and $b(z)\in z\mZ[[z]]$ such that the associated $\Phi(z)$ satisfies the following two equations:
	\begin{align*}
		&\Phi(0)=1,\\
		&\Phi(\alpha_1+\alpha_2)=-Dz.\tag{4.3}
	\end{align*}
	\end {lem}
	
	{\bf Proof.} Taking $K=Z$. See that $[L(0)]=1$ and $$[L(\alpha_1+\alpha_2)]=\sum_{\alpha\in\Delta}e^\alpha+2.$$
	Since $\lambda(lz^n)=(1-z^n)^l$, then
	\begin{align*}
		&\lambda\left(a(z)[L(\alpha_1+\alpha_2)]+b(z)\right)\\
		=&\lambda\left(\sum_{n\geq 1}a_n\left(\sum_{\alpha\in\Delta}e^\alpha+2\right)z^n+\sum_{n\geq 1}b_nz^n\right)\\
		=&\prod_{n\geq 1}\left(\prod_{\alpha\in\Delta}\lambda(a_ne^\alpha z^n)\right)\cdot\lambda\left((2a_n+b_n)z^n\right)\\
		=&\Phi(z),
	\end{align*}
	by the definition of $\Phi(z)$ in $(4.2)$.  $\Box$

	We give
	\begin {conj}\label{con:3} Let $a(z)=\sum_{n\geq 1}a_nz^n$, $b(z)=\sum_{n\geq 1}b_nz^n\in z\mZ[[z]]$ be two power series defined in Lemma \ref{lem:4.4}.
	Then, $a(z)$ and $b(z)$ are given from the dimensions of homogenous components of $A(D)$ and $\angAD$.
	\end {conj}
	
	\begin{thm}\label{thm:4.2} Conjecture 2 implies Conjecture 3.
	\end{thm}
	{\bf Proof.} Assuming the Conjecture 2 holds, we have that
	$$[A(D)]={\mathcal A}(D),\ \ [\angAD]={\mathcal B}(D).$$
	Restricting on $Z$, 
	the two equations $(4.3)$ in Lemma \ref{lem:4.4} are satisfied by $a(z)={\mathcal A}(D)|_Z,  b(z)={\mathcal B}(D)|_{Z}$. $\Box$

\section{Evidences for the conjectures}

In this section, we will give some facts to support previous conjectures.

\subsection{The case $D=1$}

For $D=1$, $A(1)$ becomes the free commutative associative algebra $xk[x]\simeq t^{-1}k[t^{-1}]$ and the homology module $H_k(\AG(A(1)))$ can be computed applying the theory in \cite{GL} by imbedding $\AG(t^{-1}k[t^{-1}])$ into the affine Kac-Moody algebra $\widehat{sl_3}$ as its $F$-parabolic subalgebra ${\frak u}^-$.
Using Galand-Lepowski Theorem in \cite{GL}, we have that
\begin{thm}\label{thm:5.1}  As $sl_3(k)$-modules, one has that
	\begin{align*}&H_0(\AG(A(1)))\simeq L(0),\ \ H_1(\AG(A(1)))\simeq L(\alpha_1+\alpha_2),\\
		&H_2(\AG(A(1)))\simeq L(\alpha_1+2\alpha_2)\oplus L(2\alpha_1+\alpha_2),\\
		&H_3(\AG(A(1)))\simeq L(2\alpha_1+3\alpha_2)\oplus L(3\alpha_1+2\alpha_2).\end{align*}
	Furthermore, the conjecture 1 holds when $D=1$. $\Box$
\end{thm}
We postpone the proof of this theorem with some preliminaries in Appendix A.

By the above Theorem and Lemma \ref{lem:3.4}, 
$\mathcal{B}A(1)=\text{Inner}\ A(1)$. Since $A(1)$ is commutative,  $\text{Inner}(A(1))=0$. Then, the dimensions $a_n(1)=1$ and $b_n(1)=0$ for all $n\geq 1$, 

\subsection{The case $D=2$}

When $D\geq 2$, the explicit structure of $H_r(\AGAD)$ is unknown yet. We attack the Conjecture 3 with some concrete computation and show that the dimensions $a_n(D)$ and $b_n(D)$ give the correct number in the Conjecture 3 up to $n\leq 3$.

By Artin's Theorem, the alternative algebras generated by two elements are associative. Hence, $A(2)=T(2)$, the free non-unital associative algebra(tensor algebra) generated by two variables. Hence, $a_n(2)=2^n, n\geq1$. Meanwhile, it is easy to see that
$$b_1(2)=0,b_2(2)=1,b_3(2)=4$$

Then,
\begin{align*}
	&\Phi^{(2)}(z)\equiv\left(\prod_{\alpha\in\Delta}(1-e^\alpha z)^2\right)(1-z)^4\left(\prod_{\alpha\in\Delta}(1-e^\alpha z^2)^4\right)(1-z^2)^9\\
	&\left(\prod_{\alpha\in\Delta}(1-e^\alpha z^3)^8\right)(1-z^2)^{20}.\ (\bmod\ z^4)
\end{align*}

We have that
\begin{thm}\label{thm:5.2} When $D=2$, in the ring $\mZ[t_1^{\pm1},t_2^{\pm1}]^{S_3}[[z]]/(z^4)$, one has that
	\begin{align*} 
		&\Phi^{(2)}\equiv\\
		&1-2z\text{ch }L(\alpha_1+\alpha_2)+z^2\left(\text{ch }L(2\alpha_1+2\alpha_2)+3(\text{ch }L(2\alpha_1+\alpha_2)+\text{ch }L(\alpha_1+2\alpha_2))\right)\\
		&-2z^3((\text{ch } L(2\alpha_1+3\alpha_2)+\text{ch } L(3\alpha_1+2\alpha_2))-(\text{ch } L(2\alpha_1+\alpha_2)+\text{ch } L(\alpha_1+2\alpha_2))\\
		&(\bmod\ z^4)
	\end{align*}
	Hence,
	$$\Phi^{(2)}(0)\ \equiv1(\bmod\ z^4),\ \ \Phi^{(2)}(\alpha_1+\alpha_2)=-2z\ (\bmod\ z^4)$$
	hold.  $\Box$
\end{thm}
The details of the computation can be seen in Appendix B.

\subsection{The cases $D\geq3$}

For general $D\geq 3$, we have 
\begin{pro}\label{pro:5.1}  When $D\geq 3$,
	
	(1) $\dim A(D)_1=D, \dim\angAD_1=0$,
	
	(2) $\dim A(D)_2=D^2, \dim\angAD_2={D\choose 2}$,
	
	(3) $\dim A(D)_3=D^3+{D\choose 3}, \dim\angAD_3=D^3-D^2-2{D\choose 3}$.
\end{pro}
{\bf Proof.} Let $\{x_1,\cdots,x_D\}$ be the generators of $A(D)$ with degree $1$. Consider the gradings
$$A(D)=\oplus_{n>0}A(D)_n,\ \ \ \angAD=\oplus_{n>0}\angAD_n.$$

Then, $A(D)_1=\text{span}_k\{x_1,\cdots,x_D\}$ and $\dim A(D)_1=D$.
By Artin's Theorem, any alternative algebra generated by less than three elements is associative. Hence,
$\dim A(D)_2=\dim T(D)_2=D^2$.

For three generators $x,y,z$ in $\{x_1,\cdots,x_D\}$, if $\#\{x,y,z\}<3$, the multiplication is associative and the dimension of the subspace of $A(D)_3$ spanned by them is
$$D+6{D\choose 2}=3D^2-2D.$$
By alternative laws, if $x,y,z$ are distinct generators,  the $7$ elements $$(xy)z,(xz)y,(yz)x,(yx)z,(zx)y,(zy)x,x(yz)$$ are linearly independent and span the subspace generated by $x,y,z$  in $A(D)_3$. Hence, the total dimension of $A(D)_3$ is 
$$3D^2-2D+7{D\choose 3}=D^3+{D\choose 3}.$$  

To compute $\dim\angAD_n$, we see that $$\angAD_n=(\oplus_{n_1+n_2=n,n_1\geq n_2\geq 1}A(D)_{n_1}\wedge A(D)_{n_2})/I'_n,$$ where $I'_n$ is the $k$-space spanned by homogenous  $xy\otimes z+yz\otimes x+zx\otimes y$ with degree $n$. Hence, $\dim\angAD_1=0$. And, $$\angAD_2\simeq \wedge^2 A(D)_1$$
which yields $\dim\angAD_2={D\choose 2}$.
Finally, $$\angAD_3\simeq A(D)_2\otimes A(D)_1/I'_3,$$
where $I'_3$ is spanned by $xy\otimes z+yz\otimes x+zx\otimes y$ for $x,y,z\in\{x_1,\cdots,x_D\}$.
One has that 
$$\dim I'_3=D+2{D\choose 2}+2{D\choose 3}=D^2+2{D\choose 3}$$
and 
$$\dim\angAD_3=D^3-D^2-2{D\choose3}.\ \ \ \Box$$ 

Note that the formulae are also valid when $D=2$ since ${2\choose 3}=0$. 

By Proposition \ref{pro:5.1}, we get
\begin{thm}\label{thm:5.3} When $D\geq 3$, in the ring $\mZ[t_1^{\pm1},t_2^{\pm1}]^{S_3}[[z]]/(z^4)$,
	\begin{align*} 
		&\Phi^{(D)}\equiv\\
		&1-Dz\text{ch }L(\alpha_1+\alpha_2)\\
		&+z^2\left({D\choose 2}\text{ch }L(2\alpha_1+2\alpha_2)+{D+1\choose 2}(\text{ch }L(2\alpha_1+\alpha_2)+\text{ch }L(\alpha_1+2\alpha_2))\right)\\
		&-z^3({D\choose 3}\text{ch } L(3\alpha_1+3\alpha_2)+\frac{D^3-D}{3}\\&((\text{ch } L(3\alpha_1+2\alpha_2)+\text{ch } L(2\alpha_1+3\alpha_2))-(\text{ch } L(2\alpha_1+\alpha_2)+\text{ch } L(\alpha_1+2\alpha_2))))\\
		&(\bmod\ z^4).
	\end{align*}
	We also have that
	$$\Phi^{(D)}(0)\ \equiv1(\bmod\ z^4),\ \ \Phi^{(D)}(\alpha_1+\alpha_2)=-Dz\ (\bmod\ z^4).\ \ \Box$$  \end{thm}

More details are in Appendix C.

In particular, when $D=3$, we have 
that $$a_1=3, a_2=9, a_3=28, b_1=0, b_2=3, b_3=16.$$
We have that
\begin{align*}
	&\Phi^{(3)}\equiv\\
	&1-3z\text{ch }L(\alpha_1+\alpha_2)\\
	&+3z^2\left(\text{ch  }L(2\alpha_1+2\alpha_2)+2 (\text{ch }L(2\alpha_1+\alpha_2)+\text{ch }L(\alpha_1+2\alpha_2))\right)\\
	&-z^3(\text{ch } L(3\alpha_1+3\alpha_2)+8(\text{ch } L(3\alpha_1+2\alpha_2)+\text{ch } L(2\alpha_1+3\alpha_2)\\
	&-(\text{ch } L(2\alpha_1+\alpha_2)+\text{ch } L(\alpha_1+2\alpha_2)))) (\bmod\ z^4).
\end{align*}

{\bf Remark.} In above theorem, we see that it is the free alternative algebra $A(D)$ but not the tensor algebra $T(D)$ offers the right value of $a_3(D)$.

\section{The cyclicity of the functor $A$}

Let $V$ be a vector space over $k$. $A(V)$ is the free alternative algebra generated by $V$. 
It is well known that $A$ can be seen as an analytic functor from the category of vector spaces over $k$ 
to the category of alternative algebras over $k$, which has a decomposition
$$A=\oplus_{n\geq 1}A_n, V\mapsto A(V)=\oplus_{n\geq 1}A_n(V),$$
and $A(D)$ is the previously defined free alternative algebra over $D$ generators. The suspensions of analytic functors were introduced in \cite{KM} and \cite{Ma}. In this section, we will prove that the functor $A$ is cyclic, namely that it is a suspension of some analytic functor.

\subsection{The cyclic functors and cyclic triples}

For an analytic functor $F$, 
$\Sigma F$ is the suspension of $F$. The following results are elementary.

\begin{pro}\label{pro:6.1}If $F=\oplus_{n\geq 0}F_n,G=\oplus_{n\geq 0}G_n$ are two analytic functors, then
	
	(1) $\Sigma F=0$ if and only if $F$ is a constant functor,
	
	(2) $\Sigma(\text{Id})$ is the constant functor with image $k$,
	
	(3) $\Sigma(F\oplus G)\simeq_N\Sigma F\oplus\Sigma G$,
	
	(4) $\Sigma(F\otimes G)\simeq_N\Sigma F\otimes G\oplus F\otimes\Sigma G$,
	
	(5) $\Sigma(F\circ G)\simeq_N\Sigma F\circ G\otimes\Sigma G$,
	
	(6) Let $\Theta:F\rightarrow G$ be a homogenous natural transform. If $F(0)\simeq G(0)$ and $\Sigma\Theta:\Sigma F\rightarrow\Sigma G$ is an isomorphism, then $\Theta$ is also an isomorphism.
\end{pro}
{\bf Remark: }The properties of suspensions are very much like those of the derivations of smooth functions.

\begin{defn}\label{defn:6.1}Let $F$ be an analytic functor. If there is another analytic functor $G$ such that $F=\Sigma G$, then $F$ is called a cyclic functor. 
\end{defn}
Let $F,G$ be two analytic functors and let $\Theta:F\rightarrow G$ be a homogenous natural transformation. Then, $\ker\Theta$ is also an analytic functor.
Furthermore, by chasing on the commutative diagram,
there are natural isomorphisms
$$\Sigma\ker\Theta\simeq_N\ker\Sigma\Theta, \ \ \Sigma\text{coker}\Theta\simeq_N\text{coker}\Sigma\Theta,$$
i.e., the suspension commutes with the kernel and the cokernel.

Let $F,G$ be two analytic functors and let $\Theta:F\otimes\text{Id}\rightarrow G$ be a homogenous natural transform, where $\text{Id}$ is the identity functor. By (4) of Proposition \ref{pro:6.1},
$$\Sigma(F\otimes\text{Id})\simeq_N\Sigma F\otimes\text{Id}\oplus F\otimes C$$ 
where $C$ is the constant functor with image $k$.
Taking value on $k^D$,
$$\Sigma(F\otimes\text{Id})(D)\simeq\Sigma F(D)\otimes k^D\oplus F(D).$$
Hence, $\Sigma F(D)\otimes k^D$ can be viewed as a summand of $\Sigma(F\otimes\text{Id})(D)$.

\begin{defn}\label{defn:6.2}The triple $(F,G,\Theta)$ is called a cyclic triple if the induced map obtained from the composition
	$$\Sigma F(D)\otimes k^D\hookrightarrow\Sigma(F\otimes\text{Id})(D)\stackrel{(\Sigma\Theta)_{k^D}}{\rightarrow}\Sigma G(D)$$
	is an isomorphism, for any $D\geq 0$.
\end{defn}

\subsection{Cyclic structures on $A$ and ${\mathcal B}A$}
For a $k$-space $V$, recall that ${\mathcal B}A(V)$ is functorial on $V$. The natural map $A(V)\otimes V\rightarrow{\mathcal B}A(V), a\otimes v\mapsto\langle a,v\rangle$, defined for $V\in\text{Vect}_k$  is indeed a natural transformation $$\Theta_A:A\otimes\text{Id}\rightarrow{\mathcal B}A.$$

$(A,{\mathcal B}A,\Theta_A)$ is a cyclic triple, which is the main result of this section. That is,
\begin{thm}\label{thm:6.1} We have an isomorphism
	$$\Sigma A(D)\otimes k^D\rightarrow\Sigma\BAD,$$
	i.e., the triple $(A,{\mathcal B}A,\Theta_A)$ is cyclic.
\end{thm}
Moreover, by Lemma 16 in \cite{KM},
\begin{cor}\label{cor:6.1} The functor $A$ is a cyclic functor. $\Box$ 
\end{cor}

\

We need some preparations to prove the above theorem.

First, recall that

{\bf Schreier-Kurosh-Cohn Theorem } Let $F$ be a free
Lie algebra over $k$ and $M$ a free $F$-module. Then every submodule $N$ of $M$ is also free.

Let $D\geq 1$. Set $F_1=sl_3(k)\otimes k^D$. $F_1$ is an $sl_3(k)$-module isomorphic to $D$ copies of the adjoint module $L(\alpha_1+\alpha_2)$. Let $F$ be the free Lie algebra generated by $F_1$. Then, the actions of $sl_3(k)$ on $F_1$ are naturally extended on $F$ such that
they are derivations of $F$.

Let $M$ be an $F$-module and a weight $sl_3(k)$-module simultaneously. And, assume that these actions are compatible. These modules are objects  of the category ${\mathcal M}(F,sl_3(k))$, whose homomorphisms are morphisms of both $F$-modules and $sl_3(k)$-modules.

Since $M$ is a weight $sl_3(k)$-module,
$$M=M^{sl_3}\oplus sl_3(k).M.$$
We have that
$$H_0(sl_3,M)=M/sl_3(k).M\simeq M^{sl_3}=H^0(sl_3,M).$$

Then,
$$M\otimes F\simeq F\otimes M\rightarrow M, g\otimes m\mapsto g.m$$
is a morphism of $sl_3(k)$-modules.
Consider the $sl_3(k)$-submodules $F_1\subseteq F$ and $M^{\text{ad}}\subseteq M$. The morphism of $sl_3(k)$-modules
$$M^{\text{ad}}\otimes F_1\hookrightarrow M\otimes F\rightarrow M$$
induces a $k$-linear morphism 
$$\mu_M:H_0(sl_3,M^{\text{ad}}\otimes F_1)\rightarrow H_0(sl_3,M).$$  

One has that
\begin{lem}\label{lem:6.1}Let $\alpha$ be a dominant integral weight in the root lattice of type $A_2$. Then, $[L(\alpha_1+\alpha_2)\otimes L(\alpha):L(0)]\neq 0$ if and only if  $\alpha=\alpha_1+\alpha_2$.  
\end{lem}
{\bf Proof. }For $\alpha=\alpha_1+\alpha_2$, recall that
\begin{equation*}L(\alpha)^{\otimes2}\simeq L(2\alpha)\oplus L(\alpha+\alpha_1)\oplus L(\alpha+\alpha_2)\oplus 2L(\alpha)\oplus L(0).\tag{6.1}\end{equation*}
Hence, $[L(\alpha)^{\otimes2}:L(0)]=[k]\neq0$.

On the other hand, by Exercise 12 of Section 24 in \cite{Hu}, the only possible dominant integral weights $\beta$ for which $L(\beta)$ can occur as a summand of $L(\alpha_1+\alpha_2)\otimes L(\alpha)$ are those in the set $\Pi(\alpha_1+\alpha_2)-\alpha$, where
$\Pi(\alpha_1+\alpha_2)=\Delta\cup\{0\}$ is the set of weights in $L(\alpha_1+\alpha_2)$. Hence, if $0\in\Pi(\alpha_1+\alpha_2)-\alpha$,
we must have that $$\alpha\in\Pi(\alpha_1+\alpha_2)\cap\Lambda^+=\{0,\alpha_1+\alpha_2\}.$$ But when $\alpha=0$, $[L(\alpha_1+\alpha_2):L(0)]=0$. The only possibility is $\alpha=\alpha_1+\alpha_2$. $\Box$

\begin{thm}\label{thm:6.2} Let
	$$0\rightarrow Y\rightarrow X\rightarrow M\rightarrow0$$
	be a short exact sequence in ${\mathcal M}(F,sl_3(k))$. Assume that
	
	1. $X$ is a free $F$-module generated by $sl_3(k).X$,
	
	2. $Y$ is generated by $sl_3(k).Y$,
	
	Then, the map $\mu_M$ is an isomorphism of $k$-spaces.
\end{thm}
{\bf Proof. }Since $F$ is freely generated by $F_1$, the universal enveloping algebra $U(F)$ is isomorphic to $T(F_1)$ by Yoneda Lemma. Moreover, by PBW Theorem, its argumentation ideal $F.U(F)\simeq F_1T(F_1)$.

Set $S=sl_3(k).X$, which is an $sl_3(k)$-module. Since $X$ is a free $F$-module generated by $S$, we have that $$F.X\simeq FU(F)\otimes S, \ \ X\simeq F.X\oplus S$$ and $$X\simeq U(F)\otimes S\simeq S\otimes T(F_1)$$ as $sl_3(k)$-modules.
Hence, \begin{align*}
	&X\otimes F_1\simeq S\otimes T(F_1)\otimes F_1 \simeq S\otimes F_1T(F_1)\\
	&\simeq S\otimes F.U(F)\simeq F.X.
\end{align*}
Particularly, the map
\begin{equation*}
H_0(sl_3,X\otimes F_1)\rightarrow H_0(sl_3,F.X).\tag{6.2}
\end{equation*}
is an isomorphism.

Since $F_1$ is an $sl_3(k)$-module of adjoint type,
$F_1=F_1^{\text{ad}}$. By Lemma \ref{lem:6.1}, when $\alpha\neq\alpha_1+\alpha_2$,
$$[L(\alpha)\otimes L(\alpha_1+\alpha_2):L(0)]=0,$$ we have that
\begin{equation*}H_0(sl_3,X\otimes F_1)=(X\otimes F_1)^{sl_3}=(X^{\text{ad}}\otimes F_1)^{sl_3}\simeq H_0(sl_3,X^{\text{ad}}\otimes F_1).\tag{6.3}\end{equation*}

Meanwhile, note that $X=F.X\oplus S$. Then,
$X/S\simeq F.X$ as $sl_3(k)$-modules.
Since $H_0$ is right exact, the short exact sequence of $sl_3(k)$-modules $$0\rightarrow S\rightarrow X\rightarrow X/S\rightarrow0$$ 
gives an exact sequence of $k$-spaces
$$\cdots\rightarrow H_0(sl_3, S)\rightarrow H_0(sl_3, X)\rightarrow H_0(sl_3, X/S)\rightarrow0.$$ 
The fact $H_0(sl_3,S)\simeq S^{sl_3}=(sl_3(k).X)^{sl_3}=0$ yields \begin{equation*}H_0(sl_3, X)\simeq H_0(sl_3, X/S)\simeq H_0(sl_3,F.X).\tag{6.4}\end{equation*}
From $(6.2), (6.3)$ and $(6.4)$,
$$H_0(sl_3,X^{\text{ad}}\otimes F_1)\simeq H_0(sl_3,X\otimes F_1)\simeq H_0(sl_3,F.X)\simeq H_0(sl_3,X)$$ and $$\mu_X:H_0(sl_3,X^{\text{ad}}\otimes F_1)\rightarrow H_0(sl_3,X)$$
is an isomorphism.  

Finally, by Schreier-Kurosh-Cohn Theorem, $Y$ is also free, and therefore $\mu_Y$ is an isomorphism.  We have the following commutative diagram
$$\xymatrix{
	H_0(sl_2, Y^{\text{ad}}\otimes F_1)\ar[d]^{\mu_X}\ar[r]
	& H_0(sl_2, X^{\text{ad}}\otimes F_1) \ar[d]^{\mu_Y}\ar[r]& H_0(sl_2, M^{\text{ad}}\otimes F_1)\ar[d]^{\mu_M}\ar[r]&0 \\
	H_0(sl_2, Y) \ar[r]& H_0(sl_2, X) \ar[r]& H_0(sl_2, M)\ar[r]&0}.
$$

By the snake lemma, $\mu_M$ is also an isomorphism. $\Box$

\

Let $M$ be an $F$-module with admissible $sl_3(k)$-actions. Set 
$M^{\text{higher}}$
is the $F$-submodule of $M$ generated by the $L(\alpha)$-components of $sl_3(k)$ for $\alpha\neq0,\alpha_1+\alpha_2$. Then, 
$$M=M^{\text{higher}}+(M^{\text{ad}}\oplus M^{sl_3}).$$
And, $M/M^{\text{higher}}$ is an $F$-module in the category $\mathfrak R$.

Note that $M^{\text{ad}}$ and $M^{sl_3}$ may not be $F$-submodules and the sum between $M^{\text{higher}}$ and $M^{\text{ad}}\oplus M^{sl_3}$ is not necessarily direct.

\

Since the Lie algebra $\AGAD$ is generated by $sl_3(k)\otimes A_1(D)\simeq sl_3(k)\otimes k^D$, there is a surjective homomorphism  of Lie algebras
$$\Psi:F\rightarrow\AGAD.$$
Furthermore, the $\AGAD$-modules with admissible $sl_3(k)$-actions are also $F$-modules in the category ${\mathcal M}(F,sl_3(k))$.

\begin{thm}\label{thm:6.3}Let $M$ be a free $\AGAD$-module in the category ${\mathcal M}_{\mathfrak R}(\AGAD)$ generated by one copy of $L(\alpha_1+\alpha_2)$.
	Then the map 
	$$\mu_M:H_0(sl_3,M^{\text{ad}}\otimes F_1)\rightarrow H_0(sl_3,M)$$ is an isomorphism.
\end{thm}
{\bf Proof.} Let $X$ be the free $F$-module generated by $L(\alpha_1+\alpha_2)$.  Set 
$X^{\text{higher}}$
is the submodule of $M$ generated by its $L(\alpha)$-components other than $\alpha=0,\alpha_1+\alpha_2$ and $P=X/X^{\text{higher}}$.
The quotient of $F$-modules is denoted by
$\pi:X\rightarrow P$. $P$ is in the category ${\mathcal M}_{\mathfrak R}(F)$ and it is free in ${\mathcal M}_{\mathfrak{R}}(F)$. Hence, there is also a morphism of $F$-modules $\sigma:P\rightarrow M$.

Let
$K$ be the sum of $L(\alpha)$-components for $\alpha\neq0,\alpha_1+\alpha_2$ in $F$. Since $\AGAD$ is free in the category $\text{Lie}_{\mathfrak R}$, $\AGAD=F/R$ where $R$ is the ideal of
$F$ generated by $K$. Therefore, $\ker\sigma$ is the $F$-submodule of $P$ generated by $K.P$.

Note that $P$ is an object of $\mathfrak{R}$, $P=P^{sl_3}\oplus P^{\text{ad}}$. Since $K.X^{sl_3}\subseteq X^{\text{higher}}$, one has that $K.P^{sl_3}=0$. Meanwhile, since
$$[L(\alpha)\otimes L(\alpha_1+\alpha_2):L(0)]=0,$$
for $\alpha\neq\alpha_1+\alpha_2$ by Lemma \ref{lem:6.1},
we have that $K.X^{\text{ad}}\subseteq X^{\text{higher}}+X^{\text{ad}}$ and $K.P^{\text{ad}}\subseteq P^{\text{ad}}$. Hence, $$K.P=K.(P^{sl_3}\oplus P^{\text{ad}})\subseteq P^{\text{ad}},$$
which means that $\ker\sigma$ is generated by its $L(\alpha_1+\alpha_2)$-components.

Consider the sequence in ${\mathcal M}(F,sl_3(k))$
$$X\stackrel{\pi}{\rightarrow}P\stackrel{\sigma}{\rightarrow}M.$$
Set $Y=\ker(\sigma\circ\pi)$. It follows from the descriptions of $\ker\pi$ and $\ker\sigma$ that $Y$
is generated by its $L(\alpha)$-components for $\alpha\neq 0$. Thus, the short sequence
$$0\rightarrow Y\rightarrow X\stackrel{\sigma\circ\pi}{\rightarrow }M\rightarrow0$$
satisfies the hypotheses of Theorem \ref{thm:6.2}. It follows that $\mu_M$ is an isomorphism. $\Box$

\

{\bf The proof of Theorem \ref{thm:6.1}.}
For $D\geq1$, let $A(1+D)$ be the free alternative algebra generated by $x_0,x_1,\cdots,x_D$ and $\AGADo$ its $\AG$ Lie algebra. View $\AGAD$ as the subalgebra of $\AGADo$ corresponding to $x_1,\cdots,x_D$. Since $$\AGADo={\mathcal B}A(1+D)\oplus A(1+D)\otimes sl_3(k),$$
we have 
$$\Sigma\AGAD=\Sigma{\mathcal B}A(D)\oplus\Sigma A(D)\otimes sl_3(k).$$
Moreover, $\Sigma\AGAD$ is an $\AGAD$-module in the category ${\mathcal M}_{\mathfrak{R}}(\AGAD)$. We denote it by $M$ with $$ M^{sl_3}=\Sigma{\mathcal B}A(D), \ \ M^{\text{ad}}=\Sigma A(D)\otimes L(\alpha_1+\alpha_2).$$

Let $N$ be any $\AGAD$-module in category ${\mathcal M}_{\mathfrak{R}}(\AGAD)$ generated by one copy  of the adjoint module $L(\alpha_1+\alpha_2)=\{x^{(N)}\ |\ x\in sl_3(k)\}$.
Let $\frak g$ be the Lie algebra $\AGAD\ltimes N$. Then, ${\frak g}$ is also an object in $\LieR$.

Denote the linear span of $\{x_0,x_1,\cdots,x_D\}$ over $k$ by $k^{1+D}$. 
Let $\phi$ be the $k$-linear morphism
$$\phi:k^{1+D}\rightarrow{\frak g}$$
such that $\phi|_{k^D}$ is composition of the embeddings
$$k^D\hookrightarrow A(D)\hookrightarrow \mathcal{BM}(\AGAD)\hookrightarrow\mathcal{BM}({\frak g})$$
and $\phi|_{kx_0}$ is the natural embedding $x_0\mapsto e_{13}(x_0)\mapsto e^{(N)}_{13}\in N$.

By Theorem \ref{thm:2.3}, $\AGADo$ is free in the category $\LieR$. Therefore, $\phi$ can extend to a Lie algebra morphism
$$\Phi:\AGADo\rightarrow{\frak g},$$
which is compatible with the $sl_3(k)$-actions.

Obviously, $\Phi$ sends $M$ to $N$, which means that $M$ is the free object in 
${\mathcal M}_{\mathfrak{R}}(\AGAD)$ generated by $L(\alpha_1+\alpha_2)$. 

Meanwhile, since $F_1=k^D\otimes L(\alpha_1+\alpha_2)$, we have that 
$$M^{ad}\otimes F_1\simeq\Sigma A(D)\otimes k^D\otimes(L(\alpha_1+\alpha_2)^{\otimes2}).$$
Also by the decomposition $(6.1)$,
$$(M^{ad}\otimes F_1)^{sl_3}\simeq\Sigma A(D)\otimes k^D.$$

By Theorem \ref{thm:6.3}, $\mu_M:H_0(sl_3,M^{\text{ad}}\otimes F_1)\rightarrow H_0(sl_3,M)$ is an isomorphism, which amounts to the fact that the bottom horizontal arrow in the following diagram is an isomorphism
$$\xymatrix{
	H_0(sl_3, M^{\text{ad}}\otimes F_1)\ar[d]^{=}\ar[r]^{\mu_M}
	& H_0(sl_3,M) \ar[d]^{=} \\
	(M^{ad}\otimes F_1)^{sl_3}\ar[d]\ar[r]
	&  M^{sl_3} \ar[d]^{=} \\
	\Sigma A(D)\otimes k^D\ar[r]& \Sigma{\mathcal B}A(D) }.
$$

We obtain the isomorphism
$\Sigma A(D)\otimes k^D\rightarrow\Sigma{\mathcal B}A(D)$ finally
and the triple $(A,{\mathcal B}A,\Theta_A)$ is cyclic. $\Box$

\subsection{An application of the cyclicity}

Here, we give an application of Theorem \ref{thm:6.1}.
For $r\geq 0$, let $\Lambda^r$ be the analytic functor of $r$-power wedge product. It is easy to see that $\Sigma\Lambda^r=\Lambda^{r-1}$. For $\Sigma H_*(\AGAD)$, one has that
\begin{lem}\label{lem:6.3} $\Sigma H_*(\AGAD)$ can be computed by the complex $$(\Lambda^*\AGAD\otimes\Sigma\AGAD)[1],$$ where $[1]$ is the translation functor on complexes.
\end{lem}
{\bf Proof. }Since the suspension $\Sigma$ commutes with the kernel and the cokernel, it also commutes with all $H_r$. Then, $\Sigma H_*(\AGAD)$ can be computed by the complex $\Sigma(\Lambda^*\AGAD)$

By $(5)$ of Proposition \ref{pro:6.1}, $$\Sigma(\Lambda^r\AGAD)\simeq\Lambda^{r-1}(\AGAD)\otimes\Sigma\AGAD.$$
We obtain the result. $\Box$

The following conjecture concerns only the invariant part of $H_*(\AGAD)$ with respect to the actions of $sl_3(k)$.
\begin {conj}\label{con:4}As $sl_3(k)$-modules,
$$H_r(\AGAD)^{sl_3}=0$$
for any $r\geq 1$.
\end {conj}
Obviously, since $H_1(\AGAD)\simeq L(\alpha_1+\alpha_2)\otimes k^D$, $H_1(\AGAD)^{sl_3}=0$. Conjecture \ref{con:4} is really the pleasant half of
Conjecture \ref{con:1}. On the other hand, analogue to Theorem 2 in \cite{KM},
\begin{thm}\label{thm:6.4} If Conjecture \ref{con:4} holds for $\AGADo$, then Conjecture \ref{con:2} holds for $\AGAD$.
\end{thm}
{\bf Proof. }The proof is similar to the proof of Theorem \ref{thm:4.1}. 
Assume that Conjecture \ref{con:4} holds for $\AGADo$.

Observe that both $H_r(\AGAD)$ and $\Sigma H_r(\AGAD)$  are sumands of $H_r(\AGADo), r\geq 0$. It follows that $$H_r(\AGAD)^{sl_3}=0, r\geq 1.$$ By Euler's characteristic formula, 
$$[\lambda([\AGAD]):L(0)]=[H_0(\AGAD):L(0)]=[k].$$ 

By Lemma \ref{lem:6.3}, we have that
$$\lambda([\AGAD]).[\Sigma\AGAD]=\sum_{r\geq 0}(-1)^r[\Sigma H_{r+1}(\AGAD)].$$
Hence, $$[\lambda([\AGAD]).[\Sigma\AGAD]:L(0)]=0.$$
Meanwhile, by Theorem \ref{thm:6.1}, 
\begin{align*}&[\Sigma\AGAD]=[\Sigma{\mathcal B}A(D)]+[\Sigma A(D)][L(\alpha_1+\alpha_2)]\\
	=&[\Sigma A(D)](D[L(0)]+[L(\alpha_1+\alpha_2)]).\end{align*}
Also by Lemma \ref{lem:6.1}, one has that
$$(D[\lambda([\AGAD]):L(0)]+[\lambda([\AGAD]):L(\alpha_1+\alpha_2)])[\Sigma A(D)]=0.$$
One can cancel $[\Sigma A(D)]$ in the left hand side of the equation and obtain that
$$[\lambda([\AGAD]):L(\alpha_1+\alpha_2)]=-D[\lambda([\AGAD]):L(0)]=-[k^D].$$

By Lemma \ref{lem:4.1}, we have that Conjecture \ref{con:2} holds. $\Box$

\section{Concluding Remarks}

The main part of this paper are two theorems and four conjectures.
\begin{thm}\label{thm:7.1}(Theorem of the adjoint pair) There are two functors, Allison-Benkart-Gao functor $\AG:\Alt\rightarrow\LieR$ and Berman-Moody functor $\BM:\LieR\rightarrow\Alt$, which form is a pair of adjoint functors, namely, there exists a natural equivalence
	$$\Hom_\LieR(\AG(A),{\mathfrak g})\simeq_N\Hom_\Alt(A,\BM({\mathfrak g}))$$
	for any $A\in\Alt$ and ${\mathfrak g}\in\LieR$.
\end{thm}

\begin{thm}\label{thm:7.2}(Theorem of the cyclicity) There is an isomorphism
	$$\Sigma A(D)\otimes k^D\rightarrow\Sigma{\mathcal B}A(D), a\otimes v\mapsto\{a,v\}.$$
	Hence, the triple $(A,{\mathcal B}A,\Theta_A)$ is cyclic and the alternative functor $A$ is a cyclic functor.
\end{thm}

The four conjectures are as follows.
\begin {conje}\label{con:7.1}(On the homology of $\AGAD$)As $sl_3(k)$-modules,
\begin{equation*}H_r(\AGAD)^{sl_3}=0 \text{\ \ and   }\ \ H_r(\AGAD)^{\text{ad}}=0\tag{7.1}\end{equation*}
for any $r\geq 2$.
\end {conje}

\begin {conje}\label{con:7.2}(On the classes of $A(D)$ and $\innAD$ in ${\mathcal R}(GL(D))$) Let $D$ be a non-negative integer and $A(D)$ the free alternative algebra with $D$ generators. In the Grothendieck ring ${\mathcal R}(GL(D))$, we have that
$$[A(D)]={\mathcal A}(D),\ \ [\angAD]={\mathcal B}(D),$$
where ${\mathcal A}(D)$ and ${\mathcal B}(D)$ are uniquely determined by two equations in ${\mathcal R}(GL(D))$:
\begin{align*}
	&[\lambda\left({\mathcal A}(D)[L(\alpha_1+\alpha_2)]+{\mathcal B}(D)\right):L(0)]=[k],\\
	&[\lambda\left({\mathcal A}(D)[L(\alpha_1+\alpha_2)]+{\mathcal B}(D)\right):L(\alpha_1+\alpha_2)]=-[k^D],\tag{7.2}
\end{align*}
where $k$ is the trivial $GL(D)$-module and $k^D$ is the natural $GL(D)$-module.
\end {conje}

\begin {conje}\label{con:7.3}(On the dimensions of $A(D)_n$ and $\innAD_n$) Let $a(z)=\sum_{n\geq 1}a_nz^n$, $b(z)=\sum_{n\geq 1}b_nz^n\in z\mZ[[z]]$ be the solutions of two equations:
\begin{align*}
	&\Phi(0)=1,\\
	&\Phi(\alpha_1+\alpha_2)=-Dz,\tag{7.3}
\end{align*}
where $\Phi$ is defined by
\begin{align*}
	&\Phi=\prod_{n\geq 1}\left(\prod_{\alpha\in\Delta}(1-e^\alpha z^n)^{a_n}\right)(1-z^n)^{2a_n+b_n}\\
	&\in\mZ[[z]][t_1^{\pm 1},t_2^{\pm 1}]^{S_3}
\end{align*}

Then, $a(z)$ and $b(z)$ are generated from the dimensions of the degree $n$ homogenous component of $A(D)$ and $\angAD$.
\end {conje}

\begin {conje}\label{con:7.4}(On the homology of $\AGAD$, weaker version)As $sl_3(k)$-modules,
\begin{equation*}H_r(\AGAD)^{sl_3}=0\tag{7.4}\end{equation*}
for all $D\geq 1$ and $r\geq 1$.
\end {conje}

We have proven that Conjecture \ref{con:1}  $\Rightarrow$ \ref{con:4} $\Rightarrow$ \ref{con:2}  $\Rightarrow$ \ref{con:3} in Theorem \ref{thm:4.1}, \ref{thm:4.2} and \ref{thm:6.4}. These results are parallel to the results of \cite{KM} and would lead to important new insight into the free alternative algebra in several generators if the conjectures are true at least.
Some evidences are given to support these conjectures in Section 5. 

When this paper is closely completed, the author read a new article \cite{DK} written by V. Dotsenko and I. Kashuba. They study systematically ``the three graces'', Lie, associative, and commutative associative algebras, in the Tits-Kantor-Koecher category, the category of $sl_2$-modules that are sums of copies of the trivial and the adjoint representation. Their results are attractive and inspire more deeply thought.

In some sense, our Allison-Benkart-Gao category $\mathfrak R$ concerned in present paper is an analogue of TKK category of \cite{DK} for $A_2$-type. Lie, associative and
commutative associative algebras are also  “the three graces” of ${\mathcal P}$-algebra in $\frak R$. Indeed, the free finitely generated Lie algebras in $\mathfrak R$ are $\AGAD$. The main conjecture(Conjecture 1) of our paper states that their truncated homologies in $\frak R$ are centrated in degrees $0$ and $1$. The corresponding questions for
 the two other graces in ${\mathfrak R}$, the free commutative associative algebra and the free associative algebra,  are also worthy and interesting.

\begin{appendices}

\section*{Appendix A}

Consider the untwisted affine Lie algebra of type $A_2$
$$\widehat{sl_3}=sl_3\otimes k[t^{\pm1}]\oplus kc\oplus kd,$$
which is the realization of the affine Kac-Moody algebra ${\frak g}(A)$ associated with the generalized Cartan matrix
$$A=\begin{pmatrix}
	2 & -1 & -1\\
	-1 & 2 & -1\\
	-1 & -1 &2\\
\end{pmatrix}.$$

Let $\Delta=\Delta_+\cup\Delta_-$ be the root system of $\widehat{sl_3}$, where $\Delta_+$(resp. $\Delta_-$) is the set of positive roots(resp. negative roots). $\Pi=\{\alpha_0,\alpha_1,\alpha_2\}$ is the set of simple roots and $\Pi^\vee=\{\alpha^\vee_0,\alpha^\vee_1,\alpha^\vee_2\}$ is the set of dual primitive roots. 
$\widehat{sl_3}$ has a triangular decomposition
$$\widehat{sl_3}=\oplus_{\alpha\in\Delta_-}(\widehat{sl_3})_\alpha\oplus{\frak h}^e\oplus\oplus_{\alpha\in\Delta_+}(\widehat{sl_3})_\alpha.$$ 
where ${\frak h}^e=k\{\alpha^\vee_0,\alpha^\vee_1,\alpha^\vee_2\}\oplus kd$ is its Cartan subalgebra.

Let $W$ be the Weyl group of $\widehat{sl_3}$ generated by
the reflections $r_{\alpha_i}, i=0,1,2$ and there is a symmetric bilinear form $(,)$ on $({\frak h^e})^*$, which is invariant under the action of $W$. For $w\in W$, define
$$\Phi_w=\Delta_+\cap w(\Delta_-).$$
Let $l(w)$ be the length of $w$. Then, $l(w)=|\Phi_w|$.
For a finite subset $\Phi$ of $\Delta$, define
$\langle\Phi\rangle=\sum_{\alpha\in\Phi}\alpha$.
We have $\langle\Phi_{r_{\alpha_i}w}\rangle=r_{\alpha_i}\langle\Phi_w\rangle+\alpha_i$ for any $w\in W$, $i=0,1,2$.

Let $\delta=\alpha_0+\alpha_1+\alpha_2\in{({\frak h}^e})^*$ be basic positive imaginary root. For any $w\in W$, $w(\delta)=\delta$. Fix  a weight $\rho\in({{\frak h}^e})^*$ such that $\rho(\alpha_i^\vee)=1$ for all $i=1,2,3$. Then, for all $w\in W$, $\langle\Phi_w\rangle=\rho-w(\rho)$(Proposition 2.5 in \cite{GL}). 

Let $S=\{\alpha_1,\alpha_2\}$ and $B$ the submatrix of $A$ associated with $S$, which is of finite type $A_2$. Then,
there is a natural injection $sl_3(k)\hookrightarrow\widehat{sl_3}$, where the subalgebra of $\widehat{sl_3}$ is generated by
$\{\alpha_i^\vee,e_{\alpha_i},f_{\alpha_i}\}_{i=1,2}$. 

Set
$$\Delta^S=\Delta\cap(\mZ\alpha_1+\mZ\alpha_2), \Delta^S_+=\Delta_+\cap\Delta^S,\Delta^S_-=\Delta_-\cap\Delta^S,$$
and the following subalgebras of $\widehat{sl_3}$
\begin{align*}
	&{\frak n}^\pm_S=\oplus_{\alpha\in\Delta^S_\pm}(\widehat{sl_3})_\alpha,\ \ \ {\frak u}^\pm=\oplus_{\alpha\in\Delta_\pm\backslash\Delta^S_\pm}(\widehat{sl_3})_\alpha\\
	&{\frak r}^e=sl_3(k)+{\frak h}^e=sl_3(k)\oplus k\alpha_0^\vee\oplus kd,\ \ \ {\frak p}^e={\frak r}^e\oplus{\frak u}^+.
\end{align*}
Then, ${\frak r}^e$ is a finite-dimensional reductive Lie algebra and we call ${\frak p}^e$ the $F$-parabolic subalgebra of $\widehat{sl_3}$ and $\widehat{sl_3}={\frak p}^e\oplus{\frak u}^-$. Set the subset of Weyl group $W$
$$W^1_S=\{w\in W\ |\ w^{-1}\Delta^S_+\subseteq\Delta_+\}.$$
In fact, when $w\in W^1_S$, $w^{-1}(\alpha_i)\in\Delta_+$ for all $\alpha_i\in S$. It yields that
\begin{align*}&(w(\rho)-\rho,\alpha_i)=(w(\rho),\alpha_i)-(\rho,\alpha_i)\\
	=&(\rho,w^{-1}(\alpha_i))-1\geq 0.\end{align*}
Hence, $w(\rho)-\rho$ is a dominant integral weight with respect to ${\frak r}^e$ and the simple highest weight module $L(w(\rho)-\rho)$ of ${\frak r}^e$ is finite-dimensional. Note that $$w(\rho)-\rho\in \text{span}_\mZ\{\alpha_0,\alpha_1,\alpha_2\}.$$ 
One has that $w(\rho)-\rho=\lambda+n\delta$ for $\lambda\in P_+$ is dominant integral and $n\leq 0$.
The multiplicities of $\delta$ indicate the degree in $t^{-1}$. We also see that for two weights $\lambda,\lambda'\in\text{span}_\mZ\{\alpha_0,\alpha_1,\alpha_2\}$, $L(\lambda)\simeq L(\lambda')$ as $sl_3(k)$-modules if and only if $\lambda\equiv\lambda'(\bmod\ \mZ\delta)$.

Obviously, 
$${\frak u}^-=\{x\in\widehat{sl_3}\ |\ [d,x]=nx \text{ for } n<0\},$$
and ${\frak u}^-=sl_3(k)\otimes(t^{-1}k[t^{-1}])$ in the realization of $\widehat{sl_3}$. The following theorem gives the structure of $H_j({\frak u}^-), j\geq 0$.
\begin{thm}\label{thm:A.1} The $j$-th homology space $H_j({\frak u}^-)$ is finite-dimensional, and it is naturally ${\frak r}^e$-module isomorphic to the direct sum
	$$\oplus_{w\in W^1_S,l(w)=j}L(w(\rho)-\rho).$$ 
\end{thm}
{\bf Proof.} Take the quasi-simple module $X$ to be the trivial $\widehat{sl_3}$-module $k$ with highest weight $0$ 
in Theorem 8.6 of \cite{GL}. $\Box$

\

{\bf The proof of Theorem \ref{thm:5.1}} 

It is easy to observe that
\begin{align*}&\{w\in W_S^1\ |\ l(w)=0\}=\{1\}, \ \ \{w\in W_S^1\ |\ l(w)=1\}=\{r_{\alpha_0}\}\\
	&\{w\in W_S^1\ |\ l(w)=2\}=\{r_{\alpha_0}r_{\alpha_1},r_{\alpha_0}r_{\alpha_2}\},\\
	&\{w\in W_S^1\ |\ l(w)=3\}=\{r_{\alpha_0}r_{\alpha_1}r_{\alpha_2},r_{\alpha_0}r_{\alpha_2}r_{\alpha_1}\}.
\end{align*}

Obviously, $H_0(\AG(A(1)))\simeq L(0)$.

Since $r_{\alpha_0}(\rho)-\rho=-\alpha_0=(\alpha_1+\alpha_2)-\delta$, $H_1(\AG(A(1)))\simeq L(\alpha_1+\alpha_2)$.

And, $r_{\alpha_0}r_{\alpha_1}(\rho)-\rho=\alpha_1+2\alpha_2-2\delta$. Hence,
$$H_2(\AG(A(1)))\simeq L(\alpha_1+2\alpha_2)\oplus L(2\alpha_1+\alpha_2).$$

Similarly,
$r_{\alpha_0}r_{\alpha_1}r_{\alpha_2}(\rho)-\rho=2\alpha_1+3\alpha_2-4\delta$. Hence,
$$H_3(\AG(A(1)))\simeq L(2\alpha_1+3\alpha_2)\oplus L(3\alpha_1+2\alpha_2).$$
We obtain the first part of this theorem. 

To show the conjecture 1 holds, it is sufficient to show that for any $w\neq w'\in W$, $$w(\rho)-\rho\not\equiv w'(\rho)-\rho\ (\bmod\ \mZ\delta).$$
Otherwise, there exits a $1\neq w\in W$ such that 
$\rho-w(\rho)\in\mZ\delta$ since $\delta$ is a fixed point for all $w\in W$. Then, $\langle\Phi_w\rangle=\rho-w(\rho)$ is a sum of imaginary roots. By Proposition 2.4 and Corollary 2.6 in \cite{GL}, we obtain that $\langle\Phi_w\rangle=0$ and $w(\rho)=\rho$, which implies $w=1$. It is a contradiction. $\Box$

\section*{Appendix B}

The following facts of simple $sl_3(k)$-modules are needed for our computation. First, the Weyl group of $sl_3(k)$ is isomorphic to $S_3$ and the characters of simple $sl_3(k)$-modules are invariant under $S_3$. For convenience, denote the orbits of weights under the Weyl group $S_3$ as follows
\begin{align*}
	&\Delta_0=\Delta=\{\alpha_1,\alpha_2,\alpha_1+\alpha_2,-\alpha_1,-\alpha_2,-(\alpha_1+\alpha_2)\},\\
	&\Delta_1=\{2\alpha_1+\alpha_2,-\alpha_1+\alpha_2,-(\alpha_1+2\alpha_2)\},\\
	&\Delta_2=\{\alpha_1+2\alpha_2,\alpha_1-\alpha_2,-(2\alpha_1+\alpha_2)\},\\
	&\Delta_3=\{3\alpha_1+2\alpha_2,3\alpha_1+\alpha_2,-\alpha_1+2\alpha_2,-\alpha_1-3\alpha_2,-2\alpha_1+\alpha_2,-2\alpha_1-3\alpha_2\},\\
	&\Delta_4=\{2\alpha_1+3\alpha_2,2\alpha_1-\alpha_2,\alpha_1-2\alpha_2,\alpha_1+3\alpha_2,-3\alpha_1-\alpha_2,-3\alpha_1-2\alpha_2\}. \tag{B.1}
\end{align*}

Next, some characters of simple $sl_3(k)$-modules can be expressed by
\begin{align*}
	&\text{ch } L(\alpha_1+\alpha_2)=\sum_{\alpha\in\Delta_0}e^{\alpha}+2,\\	
	&\text{ch } L(2\alpha_1+\alpha_2)=\sum_{\alpha\in\Delta_1}e^{\alpha}+\sum_{\alpha\in\Delta_0}e^{\alpha}+1,\\	
	&\text{ch } L(\alpha_1+2\alpha_2)=\sum_{\alpha\in\Delta_2}e^{\alpha}+\sum_{\alpha\in\Delta_0}e^{\alpha}+1,\\	
	&\text{ch } L(2\alpha_1+2\alpha_2)=\sum_{\alpha\in\Delta_0}e^{2\alpha}+\sum_{\alpha\in\Delta_1\cup\Delta_2}e^{\alpha}+2\sum_{\alpha\in\Delta_0}e^{\alpha}+3,\\	
	&\text{ch } L(3\alpha_1+2\alpha_2)=\sum_{\alpha\in\Delta_3}e^{\alpha}+\sum_{\alpha\in\Delta_0}e^{2\alpha}+\sum_{\alpha\in\Delta_2}e^{\alpha}+2\sum_{\alpha\in\Delta_1}e^{\alpha}+2\sum_{\alpha\in\Delta_0}e^{\alpha}+2,\\
	&\text{ch } L(2\alpha_1+3\alpha_2)=\sum_{\alpha\in\Delta_4}e^{\alpha}+\sum_{\alpha\in\Delta_0}e^{2\alpha}+\sum_{\alpha\in\Delta_1}e^{\alpha}+2\sum_{\alpha\in\Delta_2}e^{\alpha}+2\sum_{\alpha\in\Delta_0}e^{\alpha}+2,\\
	&\text{ch } L(3\alpha_1+3\alpha_2)=\sum_{\alpha\in\Delta_0}e^{3\alpha}+\sum_{\alpha\in\Delta_3\cup\Delta_4}e^{\alpha}+2\sum_{\alpha\in\Delta_0}e^{2\alpha}+2\sum_{\alpha\in\Delta_1\cup\Delta_2}e^{\alpha}+3\sum_{\alpha\in\Delta_0}e^{\alpha}+4.\tag{B.2}
\end{align*}

In the case $D=2$, the computation is shown as follows.
\begin{align*}
	&\prod_{\alpha\in\Delta}(1-e^\alpha z)^2=\prod_{\alpha\in\Delta^+}(1-(e^\alpha+e^{-\alpha})z+z^2)^2\\
	\equiv&\prod_{\alpha\in\Delta^+}\left(1-2(e^\alpha+e^{-\alpha})z+(e^{2\alpha}+4+e^{-2\alpha})z^2-2(e^\alpha+e^{-\alpha})z^3\right)\\
	\equiv&1-2\left(\sum_{\alpha\in\Delta_0}e^\alpha\right)z+\left(\sum_{\alpha\in\Delta_0}e^{2\alpha}+4\sum_{\alpha\in\Delta_1\cup\Delta_2}e^{\alpha}+4\sum_{\alpha\in\Delta_0}e^\alpha+12\right)z^2\\
	&-2\left(\sum_{\alpha\in\Delta_3\cup\Delta_4}e^{\alpha}+4\sum_{\alpha\in\Delta_0}e^{2\alpha}+2\sum_{\alpha\in\Delta_1\cup\Delta_2}e^{\alpha}+9\sum_{\alpha\in\Delta_0}e^{\alpha}+8\right)z^3\ \ (\bmod\ z^4).\tag{I}
\end{align*}
And,
\begin{align*}
	&\prod_{\alpha\in\Delta}(1-e^\alpha z)^2(1-z)^4\equiv (I)\cdot(1-4z+6z^2-4z^3)\\
	\equiv&1-2\left(\sum_{\alpha\in\Delta_0}e^\alpha+2\right)z+\left(\sum_{\alpha\in\Delta_0}e^{2\alpha}+4\sum_{\alpha\in\Delta_1\cup\Delta_2}e^\alpha+12\sum_{\alpha\in\Delta_0}e^{\alpha}+18\right)z^2\\
	&-(2\sum_{\alpha\in\Delta_3\cup\Delta_4}e^{\alpha}+12\sum_{\alpha\in\Delta_0}e^{2\alpha}+20\sum_{\alpha\in\Delta_1\cup\Delta_2}e^{\alpha}+46\sum_{\alpha\in\Delta_0}e^{\alpha}+68)z^3\\
	&(\bmod\ z^4)\tag{II}. 
\end{align*}

Meanwhile,
\begin{align*}
	&\left(\prod_{\alpha\in\Delta}(1-e^\alpha z^2)^4(1-z^2)^9\right)\left(\prod_{\alpha\in\Delta}(1-e^\alpha z^3)^8(1-z^3)^{20}\right)\\
	\equiv&\left(\prod_{\alpha\in\Delta}(1-4e^\alpha z^2)(1-9z^2)\right)\left(\prod_{\alpha\in\Delta}(1-8e^\alpha z^3)(1-20z^3)\right)\\
	\equiv&1-\left(4\sum_{\alpha\in\Delta_0}e^\alpha+9\right)z^2-\left(8\sum_{\alpha\in\Delta_0}e^\alpha+20\right)z^3\ (\bmod\ z^4).\tag{III}
\end{align*}
Since 
$$\left(\sum_{\alpha\in\Delta_0}e^\alpha\right)^2=\sum_{\alpha\in\Delta_0}e^{2\alpha}+2\sum_{\alpha\in\Delta_1\cup\Delta_2}e^{\alpha}+2\sum_{\alpha\in\Delta_0}e^{\alpha}+6,$$then, 
\begin{align*}
	&\Phi^{(2)}(z)\equiv(I)\cdot(III)\\
	\equiv&1-2\left(\sum_{\alpha\in\Delta_0}e^\alpha+2\right)z+\left(\sum_{\alpha\in\Delta_0}e^{2\alpha}+4\sum_{\alpha\in\Delta_1\cup\Delta_2}e^\alpha+8\sum_{\alpha\in\Delta_0}e^{\alpha}+9\right)z^2\\
	&-\left(2\sum_{\alpha\in\Delta_3\cup\Delta_4}e^{\alpha}+4\sum_{\alpha\in\Delta_0}e^{2\alpha}+4\sum_{\alpha\in\Delta_1\cup\Delta_2}e^{\alpha}+4\sum_{\alpha\in\Delta_0}e^{\alpha}+4\right)z^3\\
	&(\bmod\ z^4).
\end{align*}

By $(B.2)$, we have that
\begin{align*}
	&\Phi^{(2)}(z)\equiv\\
	&1-2z\text{ch }L(\alpha_1+\alpha_2)+z^2\left(\text{ch }L(2\alpha_1+2\alpha_2)+3(\text{ch }L(2\alpha_1+\alpha_2)+\text{ch }L(\alpha_1+2\alpha_2))\right)\\
	&-2z^3((\text{ch } L(2\alpha_1+3\alpha_2)+\text{ch } L(3\alpha_1+2\alpha_2))-(\text{ch } L(2\alpha_1+\alpha_2)+\text{ch } L(\alpha_1+2\alpha_2))\\
	&(\bmod\ z^4).
\end{align*}

\section*{Appendix C}

There are tedious calculations for the case $D\geq 3$.
Since \begin{align*}
	&(1-(e^\alpha+e^{-\alpha})z+z^2)^{a_1}\\
	\equiv&(1-a_1(e^\alpha+e^{-\alpha})z+\left({a_1\choose 2}(e^{2\alpha}+2+e^{-2\alpha})+a_1\right)z^2\\
	&-\left({a_1\choose 3}(e^{3\alpha}+e^{-3\alpha})+\frac{a_1^2(a_1-1)}{2}(e^\alpha+e^{-\alpha})\right)z^3)&(\bmod\ z^4),
\end{align*}
we have
\begin{align*}
	&\prod_{\alpha\in\Delta}(1-e^\alpha z)^{a_1}=\prod_{\alpha\in\Delta^+}(1-(e^\alpha+e^{-\alpha})z+z^2)^{a_1}\\
	\equiv&1-a_1\left(\sum_{\alpha\in\Delta_0}e^\alpha\right)z\\
	&+\left({a_1\choose 2}\sum_{\alpha\in\Delta_0}e^{2\alpha}+a_1^2\sum_{\alpha\in\Delta_1\cup\Delta_2}e^{\alpha}+{a_1}^2\sum_{\alpha\in\Delta_0}e^\alpha+3a_1^2\right)z^2\\
	&-({a_1\choose 3}\sum_{\alpha\in\Delta_0}e^{3\alpha}+a_1{a_1\choose 2}\sum_{\alpha\in\Delta_3\cup\Delta_4}e^{\alpha}+a_1^3\sum_{\alpha\in\Delta_0}e^{2\alpha}+a^2_1(a_1-1)\sum_{\alpha\in\Delta_1\cup\Delta_2}e^{\alpha}\\
	&+\frac{a_1^2(5a_1-1)}{2}\sum_{\alpha\in\Delta_0}e^{\alpha}+2a_1^3)z^3(\bmod\ z^4).\tag{I}
\end{align*}
And,
\begin{align*}
	&\prod_{\alpha\in\Delta}(1-e^\alpha z)^{a_1}(1-z)^{2a_1+b_1}\equiv (I)\cdot\left(1-(2a_1+b_1)z+{2a_1+b_1\choose2}z^2-{2a_1+b_1\choose3}z^3\right)\\
	\equiv&1-\left(a_1(\sum_{\alpha\in\Delta_0}e^\alpha+2)+b_1\right)z\\
	&+\left({a_1\choose2}\sum_{\alpha\in\Delta_0}e^{2\alpha}+a_1^2\sum_{\alpha\in\Delta_1\cup\Delta_2}e^\alpha+a_1(3a_1+b_1)\sum_{\alpha\in\Delta_0}e^{\alpha}+(3a_1^2+{2a_1+b_1\choose2})\right)z^2\\
	&-({a_1\choose 3}\sum_{\alpha\in\Delta_0}e^{3\alpha}+\frac{a_1^2(a_1-1)}{2}\sum_{\alpha\in\Delta_3\cup\Delta_4}e^{\alpha}+\frac{a_1(4a_1^2+a_1b_1-2a_1-b_1)}{2}\sum_{\alpha\in\Delta_0}e^{2\alpha}\\
	&+a_1^2(3a_1+b_1-1)\sum_{\alpha\in\Delta_1\cup\Delta_2}e^{\alpha}+\frac{a_1(13a_1^2-3a_1+8a_1b_1+b_1^2-b_1)}{2}\sum_{\alpha\in\Delta_0}e^{\alpha}\\
	&+a_1^2(8a_1+3b_1)+{{2a_1+b_1}\choose 3})z^3 (\bmod\ z^4)\tag{II}. 
\end{align*}

Meanwhile,
\begin{align*}
	&\left(\prod_{\alpha\in\Delta}(1-e^\alpha z^2)^{a_2}(1-z^2)^{2a_2+b_2}\right)\left(\prod_{\alpha\in\Delta}(1-e^\alpha z^3)^{a_3}(1-z^3)^{2a_3+b_3}\right)\\
	\equiv&\left(\prod_{\alpha\in\Delta}(1-a_2e^\alpha z^2)(1-(2a_2+b_2)z^2)\right)\left(\prod_{\alpha\in\Delta}(1-a_3e^\alpha z^3)(1-(2a_3+b_3)z^3)\right)\\
	\equiv&1-\left(a_2\sum_{\alpha\in\Delta_0}e^\alpha+2a_2+b_2\right)z^2-\left(a_3\sum_{\alpha\in\Delta_0}e^\alpha+2a_3+b_3\right)z^3\ (\bmod\ z^4).\tag{III}
\end{align*}
Then, 
\begin{align*}
	&\Phi^{(D)}(z)\equiv(II)\cdot(III)\\
	\equiv&1-\left(a_1(\sum_{\alpha\in\Delta_0}e^\alpha+2)+b_1\right)z\\
	&+({a_1\choose2}\sum_{\alpha\in\Delta_0}e^{2\alpha}+a_1^2\sum_{\alpha\in\Delta_1\cup\Delta_2}e^\alpha+(a_1(3a_1+b_1)-a_2)\sum_{\alpha\in\Delta_0}e^{\alpha}+(3a_1^2\\
	&+{2a_1+b_1\choose2}-2a_2-b_2))z^2\\
	&-({a_1\choose 3}\sum_{\alpha\in\Delta_0}e^{3\alpha}+\frac{a_1^2(a_1-1)}{2}\sum_{\alpha\in\Delta_3\cup\Delta_4}e^{\alpha}+\frac{a_1(4a_1^2+a_1b_1-2a_1-b_1-2a_2)}{2}\sum_{\alpha\in\Delta_0}e^{2\alpha}\\
	&+a_1(3a^2_1+a_1b_1-a_1-2a_2)\sum_{\alpha\in\Delta_1\cup\Delta_2}e^{\alpha}\\
	&+\left(\frac{a_1(13a_1^2-3a_1+8a_1b_1+b_1^2-b_1-8a_2-2b_2)}{2}-a_2(2a_1+b_1)+a_3\right)\sum_{\alpha\in\Delta_0}e^{\alpha}\\
	&+a_1^2(8a_1+3b_1)+{{2a_1+b_1}\choose 3}-6a_1a_2-(2a_1+b_1)(2a_2+b_2)+(2a_3+b_3))\\
	&(\bmod\ z^4).
\end{align*}

Substituting the values of $a_i,b_i$  for $i=1,2,3$ in Proposition \ref{pro:5.1}, 
\begin{align*}
	&\Phi^{(D)}(z)\equiv\\
	&1-Dz\text{ch }L(\alpha_1+\alpha_2)\\
	&+z^2\left({D\choose 2}\text{ch }L(2\alpha_1+2\alpha_2)+{D+1\choose 2}(\text{ch }L(2\alpha_1+\alpha_2)+\text{ch }L(\alpha_1+2\alpha_2))\right)\\
	&-z^3({D\choose 3}\text{ch } L(3\alpha_1+3\alpha_2)+\frac{D^3-D}{3}\\&((\text{ch } L(3\alpha_1+2\alpha_2)+\text{ch } L(2\alpha_1+3\alpha_2))-(\text{ch } L(2\alpha_1+\alpha_2)+\text{ch } L(\alpha_1+2\alpha_2))))\\
	&(\bmod\ z^4).
\end{align*}

\end{appendices}

\end{document}